%% file: main.tex
\documentclass[preprint,12pt,3p]{elsarticle}
\input{packages}
\begin{document}

\begin{frontmatter}



\title{Adjoint-based goal-oriented implicit shock tracking using full space mesh optimization}


\author[inst1]{Pranshul Thakur\corref{cor1}}
\ead{pranshul.thakur@mail.mcgill.ca}
\cortext[cor1]{Corresponding author}

\author[inst1]{Siva Nadarajah}
\ead{siva.nadarajah@mcgill.ca}

\affiliation[inst1]{organization={Department of Mechanical Engineering},%Department and Organization
            addressline={McGill University}, 
            city={Montreal},
            postcode={H3A 0C3}, 
            state={QC},
            country={Canada}}

\input{1Abstract}
\begin{comment}
\begin{graphicalabstract}
\includegraphics{grabs}
\end{graphicalabstract}

\begin{highlights}
\item First and second-order derivatives of the dual-weighted residual were computed.
\item Full-space optimizer was used with the objective function based on the dual-weighted residual.
\item Region of the shock important for the functional was tracked by the mesh.
\end{highlights}
\end{comment}


\end{frontmatter}

\input{2Introduction}

\input{3DGDiscretization}

\input{4GoalOrientedDerivatives}

\input{5FullSpaceOptimizer}
\input{6Algorithm}
\input{7NumericalTestCases}

\input{8Conclusion}
\input{9Acknowledgements}
\clearpage
\appendix
\input{10Appendix}
\clearpage
 \bibliographystyle{elsarticle-num} 
 \biboptions{sort&compress}
 \bibliography{references}





% \end{thebibliography}
\end{document}

%% file: packages.tex
\usepackage{amsmath}
\usepackage{amssymb}
\usepackage{float}
\usepackage{cancel}
\usepackage{comment}
\usepackage{caption}
\usepackage[export]{adjustbox}
\usepackage{subcaption}
\newcommand{\vect}[1]{\boldsymbol{\mathbf{#1}}}
\usepackage{algorithm}
\usepackage{algpseudocode}

\graphicspath{{figures/}} 
\usepackage{lineno}
\usepackage[numbers]{natbib}

\usepackage{url}

%% file: 1Abstract.tex
\begin{abstract}
Solutions to the governing partial differential equations obtained from a discrete numerical scheme can have significant errors, especially near shocks when the discrete representation of the solution cannot fully capture the discontinuity in the solution. A recent approach to shock tracking \cite{zahr1,mdg1} has been to implicitly align the faces of mesh elements with the shock, yielding accurate solutions on coarse meshes. In engineering applications, the solution field is often used to evaluate a scalar functional of interest, such as lift or drag over an airfoil. While functionals are sensitive to errors in the flow solution, certain regions in the domain are more important for accurate evaluation of the functional than the rest. Using this fact, we formulate a goal-oriented implicit shock tracking approach that captures a segment of the shock that is important for evaluating the functional. Shock tracking is achieved using Lagrange-Newton-Krylov-Schur (LNKS) full space optimizer, with the objective of minimizing the adjoint-weighted residual error indicator.  We also present a method to evaluate the sensitivity and the Hessian of the functional error. Using available block preconditioners for LNKS \cite{birosghattas1,doug} makes the full space approach scalable. The method is applied to test cases of two-dimensional advection and inviscid compressible flows to demonstrate functional-dependent shock tracking. Tracking the entire shock without using artificial dissipation results in the error converging at the orders of $\mathcal{O}(h^{p+1})$.
\end{abstract}

%% file: 2Introduction.tex
\section{Introduction}
High-order discretization methods, such as the discontinuous Galerkin (DG) scheme, are increasingly gaining attention because of their ability to compute accurate solutions on coarse meshes. However, such methods require the flow features to be well-resolved before observing the high-order accuracy of the scheme. Discontinuities, such as shocks, severely impact the stability of high-order schemes and one has to perform additional operations in order to obtain a stable numerical solution. A common approach to capturing discontinuous solutions using high-order schemes is to add artificial dissipation to the governing partial differential equation (PDE) \cite{Hartmann_euler, subcell_shock_capturing}. More recently, using the fact that the DG scheme allows the solution to be discontinuous across element's faces, shock tracking schemes have been developed that align the faces of mesh elements with the shock while using the discontinuous Galerkin scheme to discretize the PDE. In this regard, we mention the method of moving mesh discontinuous Galerkin with interface condition enforcement (MDG-ICE) \cite{mdg1,mdg2,hongluo_aiaa}. MDG-ICE solves the flow equation while enforcing the interface condition, where the interface condition refers to the condition of zero jump in
flux across elements' faces. Another method of shock tracking, introduced by Zahr and Persson \cite{zahr1}, solves an optimization problem to align the mesh faces with the shock. In \cite{zahr2}, a sequential quadratic programming (SQP) solver and an objective function based on the enriched residual were used to track the shock. The methods in \cite{mdg1,zahr2} were observed to provide highly accurate solutions on coarse meshes and have been extended to viscous flows \cite{corrigan_viscous} and flows in three dimensions \cite{zahr4,zahr5,zahr6}.         

The above-cited shock tracking approaches use full space optimization, in which the flow residual and the gradient of the objective are converged simultaneously. Full space optimization has been explored in the past for aerodynamic shape optimization \cite{doug, salas, hicken1, hicken2, one_shot, frank_shubin, birosghattas1, birosghattas2, gatsis_zingg}. It was observed in \cite{frank_shubin} that full space optimization is significantly more efficient than the reduced space approach. However, a robust full space approach requires globalization strategies that are different from those used in the classical reduced space approach. Further, the linear system to be solved at each step of the full space optimization can be ill-conditioned, requiring suitable preconditioners. The Lagrange-Newton-Krylov-Schur (LNKS) full space optimizer in the works of Biros and Ghattas \cite{birosghattas1, birosghattas2} uses preconditioners based on the reduced space SQP and is shown to achieve convergence nearly independent of the problem size. \cite{birosghattas2} also suggests several globalization strategies to converge the full space optimizer. LNKS has been used in quasi-one-dimensional test cases \cite{hicken1, hicken2} and in aerodynamic shape optimization \cite{doug}. Shi-Dong and Nadarajah \cite{doug} conducted an extensive study of the cost of solving the optimization problem, using different preconditioners proposed in \cite{birosghattas1}. By employing both reduced and full space optimization on problems of varying sizes, they found that the cost of full space optimization did not scale rapidly with problem size if a suitable preconditioner is available, compared to the reduced space Newton and the reduced space Quasi-Newton methods, which scaled significantly with size. \cite{doug} also used a globalization strategy involving pseudo-transient continuation to increase robustness of the full space method.   

Shock tracking algorithms summarized at the beginning of this section track the shock in the entire domain. However, we note that if the goal of a numerical computation is to evaluate a functional of interest, it is necessary to evaluate only a subset of the solution field accurately. Several studies have been conducted to determine the contribution to the functional error from the discretization errors in different regions of the computational domain. Let us mention the work of Becker and Rannacher \cite{rannacher}, Giles and Suli \cite{gilessuli} and Pierce and Giles \cite{piercegiles} for estimating the error in a functional in the context of finite element discretizations. These works express functional error in terms of residual perturbations weighted with adjoints, resulting in a form commonly known as the dual-weighted residual. The dual-weighted residual error indicator has been used to correct functionals \cite{piercegiles} and perform both goal-oriented fixed-fraction \cite{Hartmann_euler, nemec1, woopen_may, solin_fixedfraction} and anisotropic mesh adaptation \cite{inriagoaloriented, inriaviscous, dolejsigoal1, dolejsigoal2, dolejsigoal3, ringue, Yano, cezehp, fidkowski_chen, keigan} to obtain functional-dependent meshes.

In the current work, a method of goal-oriented shock tracking is developed to track shocks in regions important for the functional. Contrary to the residual-based approaches for shock tracking \cite{mdg1,zahr2}, the current approach uses the dual-weighted residual as the objective function. Moreover, we derive and use the second-order derivatives of the objective function within the full space LNKS framework, along with the available LNKS block preconditioners \cite{birosghattas1,doug} to make the algorithm scalable. The paper is organized as follows. In section \ref{DG}, we describe the DG discretization of the governing equations, with special emphasis on the simulation and control variables subsequently used as design variables for the optimization. Section \ref{goal_oriented_derivatives} then describes the goal-oriented error indicator used in this work, along with the means to compute its first and second-order derivatives. In section \ref{full_space}, we discuss the LNKS full space optimization employed in this work and outline the globalization strategies required to converge the optimizer. Finally, section \ref{test_cases} presents applications of the proposed method to several test cases to demonstrate the method's utility in tracking regions of shock that are important for evaluating the chosen functional of interest.

%% file: 3DGDiscretization.tex
\section{Discontinuous Galerkin discretization of the governing equations}
\label{DG}

This section describes the discontinuous Galerkin (DG) discretization \cite{hesthaven} of the governing equations as implemented in our in-house code, Parallel High-Order Library for PDEs (PHiLiP) \cite{doug}, using deal.II \cite{dealii} as the finite-element library. Consider the following system of $m$ steady-state nonlinear conservation laws on a domain $\Omega \subset \mathbb{R}^d$:
\begin{equation}
\label{gov_eq}
    \vect{\nabla} \cdot \vect{F(u)} = \vect{0}, \quad \quad \vect{u} : \mathbb{R}^d \to \mathbb{R}^m, \quad \vect{F(u)} : \mathbb{R}^m \to \mathbb{R}^{m \times d},
\end{equation}
%Siva with boundary conditions imposed at the boundary $\Gamma$. 
subject to appropriate boundary conditions imposed at the boundary $\Gamma$. $\vect{u(x)}$ denotes the conservative variables of this system. 

\subsection{Mathematical notation}
Let the physical domain $\Omega$ be decomposed as $\mathcal{T}_h = \{\kappa\}$, consisting of $N_k$ non-overlapping hexahedral or quadrilateral elements $\kappa$. We assume that $\kappa$ is an image of a reference element $\hat{\kappa}$ under a bijective mapping $\vect{x} = \Phi_{\kappa}(\vect{\xi}),\;\vect{x}\in\kappa,\;\vect{\xi}\in\hat{\kappa}$. Using the metric Jacobian $\vect{J} = \frac{\partial \Phi_{\kappa}}{\partial \vect{\xi}}$, the operators and integrals in the physical elements can be mapped to the reference element through the transformations:
\begin{equation}
\label{nabla_transform}
    \vect{\nabla} = \vect{J^{-T}}\vect{\nabla_\xi},
\end{equation}
\begin{equation}
\label{dx_vol}
    \int_{\kappa} u\vect{dx} = \int_{\hat{\kappa}}u\;\text{det}(\vect{J})\vect{d\xi},
\end{equation}
\begin{equation}
\label{dx_surface}
   \int_{\partial \kappa} u\vect{dx} = \int_{\partial \hat{\kappa}} u\| \vect{\text{\textbf{cof}}(J)} \vect{\hat{n}}\| \vect{d\xi},
\end{equation}
\begin{equation}
\label{normal}
    \vect{n} = \frac{ \vect{\text{\textbf{cof}}(J)} \vect{\hat{n}(\xi)}}{ \|\vect{\text{\textbf{cof}}(J)} \vect{\hat{n}(\xi)}\|},
\end{equation}
where $\vect{\text{\textbf{cof}}(J)}$ is the cofactor of the Jacobian, $\vect{n}$ and $\vect{\hat{n}}$ are the unit outward normals to $\partial \kappa$ and $\partial \hat{\kappa}$ respectively. Eq. (\ref{nabla_transform}) can be obtained through covariant transformations or by using chain rule, Eq. (\ref{dx_vol}) is the transformation of the differential volume element, Eq. (\ref{dx_surface}) indicates the transformation of differential area of the surface element and Eq. (\ref{normal}) can be obtained from covariant transformations or by taking the cross product of two vectors tangent to the face. For detailed derivations of the surface transformations, Eqs. (\ref{dx_surface}) and (\ref{normal}), the reader is referred to \cite[Appendix B.2.]{PhilipZ}. 

We denote the polynomial space of degree $q$ on the reference element $\hat{\kappa}$ as
\begin{equation*}
    \mathcal{P}^q(\hat{\kappa}) = \{ \xi_1^{i_1}...\xi_d^{i_d}, \quad\vect{\xi}\in\hat{\kappa}, \quad 0\leq i_k \leq q \quad k=1,..,d.\},
\end{equation*}
with $N_q = \text{dim}(\mathcal{P}^q(\hat{\kappa})) = (q+1)^d$. $\mathcal{P}^q(\hat{\kappa})$ can be viewed as the space of tensor product of polynomials of degree $q$ along each coordinate.

\subsection{Formulation of weak DG}
On an element $\kappa \in \mathcal{T}_h$, we take $\vect{u} \in [\mathcal{P}^p(\hat{\kappa})]^m$ and $\vect{x} \in [\mathcal{P}^q(\hat{\kappa})]^d$. They can be expressed in terms of basis functions and coefficients as:
\begin{equation}
\label{sol_poly}
    \vect{u}\vect{(\xi)} = \sum\limits_{j=1}^{N_p} \vect{\hat{u}}_j\varphi_j (\vect{\xi}), \quad \mathcal{P}^p(\hat{\kappa}) = \text{span}\{\varphi_i(\vect{\xi}) \}_{i=1}^{N_p}, 
\end{equation}
\begin{equation}
\label{mapping}
    \vect{x(\xi)} = \sum\limits_{j=1}^{N_q} \vect{\hat{x}}_j\phi_j (\vect{\xi}),\quad \mathcal{P}^q(\hat{\kappa}) = \text{span}\{\phi_i(\vect{\xi}) \}_{i=1}^{N_q},
\end{equation}
with the associated metric Jacobian
\begin{equation}
\label{jacobian}
\vect{J(\xi)} = \frac{\partial  \vect{x(\xi)}}{\partial \vect{\xi}} = \sum\limits_{j=1}^{N_q} \vect{\hat{x}}_j\frac{\partial\phi_j (\vect{\xi})}{\partial \vect{\xi}}.
\end{equation}
While $\vect{u}$ is represented using piecewise continuous polynomials in the entire domain $\Omega$, $\vect{x}$ is represented using continuous polynomials.

The discrete weak formulation of DG for (\ref{gov_eq}) on $\kappa$ can be stated as: find $\vect{u} \in [\mathcal{P}^p(\hat{\kappa})]^m$ satisfying
\begin{equation}
    \label{weak_form_physical}
    \int_{\partial \kappa}\vect{\varphi}\cdot \vect{\mathcal{H}(u^+, u^-, n^+)} \vect{dx} \\
     - \int_{\kappa} \vect{\nabla \varphi} : \vect{F(u)}\vect{dx} = 0 \quad \forall \vect{\varphi} \in [\mathcal{P}^p(\hat{\kappa})]^m,
\end{equation}
where superscripts $+$ and $-$ denote respectively the trace of the functions interior and exterior to the element's face $\partial \kappa$ and $$\vect{\nabla \varphi(x)} : \vect{F(u)} = \sum\limits_{j=1}^{m} \sum\limits_{k=1}^{d} \left(\vect{\nabla \varphi(x)}\right)_{jk}\left(\vect{F(u)})\right)_{jk}.$$
In Eq. (\ref{weak_form_physical}), $\vect{\mathcal{H}(u^+, u^-, n^+)}$ is the numerical flux at the face, satisfying $\vect{\mathcal{H}(u, u, n)} = \vect{F(u)\cdot n}$ and $\vect{\mathcal{H}(u^+, u^-, n^+)} = -\vect{\mathcal{H}(u^-, u^+, n^-)}$. For the purpose of shock tracking, $\vect{\mathcal{H}}$ also needs to satisfy the Rankine-Hugoniot jump conditions as the jump in the exact solution is non-zero across a shock \cite{mdg1, zahr2}. Using mapping transformations, Eqs. (\ref{nabla_transform}) to (\ref{normal}), along with the solution and mesh representations in Eqs. (\ref{sol_poly}) and (\ref{mapping}), the residual of the weak form in Eq. (\ref{weak_form_physical}) can be evaluated in the reference element $\hat{\kappa}$:
\begin{equation}
     \label{weak_form_final}
    \int_{\partial \hat{\kappa}}\vect{\varphi}\cdot \vect{\mathcal{H}(u^+, u^-, n^+)} \|\vect{\text{\textbf{cof}}(J)}\vect{\hat{n}}\| \vect{d\xi} 
      - \int_{\hat{\kappa}} \vect{J^{-T}}\vect{\nabla_\xi \varphi} : \vect{F(u)} \text{det}(\vect{J})\vect{d\xi} = 0 \quad \forall \vect{\varphi} \in [\mathcal{P}^p(\hat{\kappa})]^m.
\end{equation}
We represent the global vector of residuals across all the elements as $\vect{r} = \vect{r(\hat{u}, \hat{x})}$, where $\vect{\hat{u}}$ and $\vect{\hat{x}}$ denote the global vectors of coefficients of solution and mesh polynomials in Eqs. (\ref{sol_poly}) and (\ref{mapping}). A general scalar functional, $J(\vect{\hat{u}}, \vect{\hat{x}})$, can be evaluated as
\begin{equation}
\label{functional}
    \begin{split}
    J(\vect{\hat{u}}, \vect{\hat{x}}) 
    & = \int_\Omega j_\Omega(\vect{u}) \vect{dx} + \int_\Gamma j_\Gamma(\vect{u}) \vect{dx} \\
     & = \sum_{\kappa \in \mathcal{T}_h} \int_{\kappa} j_\Omega(\vect{u}) \vect{dx}
    + \int_{\partial \kappa \;\text{on} \; \Gamma} j_\Gamma(\vect{u}) \vect{dx} \\
    & = \sum_{\kappa \in \mathcal{T}_h} \int_{\hat{\kappa}} j_\Omega(\vect{u}) \text{det}(\vect{J}) \vect{d\xi}
    + \int_{\partial \hat{\kappa} \; \text{on} \; \Gamma} j_\Gamma(\vect{u}) \|\vect{\text{\textbf{cof}}(J)} \vect{\hat{n}}\| \vect{d\xi}, 
    \end{split}
\end{equation}
with $j_\Omega(\vect{u})$ and $j_\Gamma(\vect{u})$ depending on the chosen functional of interest. In order to simplify notations in the sections that follow, we denote $\vect{r(\hat{u}, \hat{x})} = \vect{r(u,x)}$ and $J(\vect{\hat{u}}, \vect{\hat{x}}) = J(\vect{u}, \vect{x})$. However, when we indicate differentiation with respect to $\vect{u}$ and $\vect{x}$, we actually refer to differentiation with respect to the coefficients $\vect{\hat{u}}$ and $\vect{\hat{x}}$ of their polynomial representations. In all test cases in this work, the surface and volume integrals arising in Eqs. (\ref{weak_form_final}) and (\ref{functional}) are evaluated using high-order Gauss-Legendre quadrature rule.

As a final note, we impose the prescribed boundary condition $\vect{u_\Gamma}$ using the numerical flux $\vect{\mathcal{H}(u_\Gamma, u_\Gamma, n)} = \vect{F(u_\Gamma)}\cdot \vect{n}$ at the face on the boundary to yield adjoint consistent residual discretization \cite{hartmann_adjoint_consistent}.

%% file: 4GoalOrientedDerivatives.tex
\section{Goal-oriented error indicator for shock tracking and its derivatives}
\label{goal_oriented_derivatives}
Information about the error in a computed functional is important to guide mesh optimization. Using the notations described in the previous section, let $\vect{u} \in [\mathcal{P}^p(\hat{\kappa})]^m$ satisfying $\vect{r(u,x)} = \vect{0}$. We take $\vect{U} \in [\mathcal{P}^{p+1}(\hat{\kappa})]^m$ to be the solution obtained by interpolating $\vect{u}$ with the interpolation operator $\vect{I_h}$, i.e. $\vect{U} = \vect{I_h}\vect{u}$. Taking the test functions for the DG discretization in $[\mathcal{P}^{p+1}(\hat{\kappa})]^m$, we obtain $\vect{R(U,x)} \neq \vect{0}$. Let $\vect{U_t} \in [\mathcal{P}^{p+1}(\hat{\kappa})]^m$ be the solution satisfying $\vect{R(U_t,x)} = \vect{0}$. The error in the functional due to the discretization error in the residual can be estimated as $\delta J = J(\vect{U_t},\vect{x})$ - $J(\vect{U},\vect{x})$. Using first-order Taylor approximations of $J$ and $\vect{R}$ about $\vect{U}$, $\vect{U_t}$ can be eliminated to get the error in the functional as \cite{venditti}
\begin{equation}
    \delta J = \vect{\psi^T}\vect{R} = \sum\limits_{i=1}^{n_{\text{dofs}}} \psi_iR_i,
\end{equation}
where $\vect{R} = \vect{R(U,x)}$, $n_{\text{dofs}}$ is the size of $\vect{R}$ and the adjoint solution, $\vect{\psi}$, is obtained from the discrete adjoint equation
\begin{equation}
\label{adjoint_eq}
    \vect{R_U^T \psi} = \vect{-J_U^T}.
\end{equation}
The linear system in Eq. (\ref{adjoint_eq}) is solved using GMRES \cite{gmres} with a tolerance of machine precision to evaluate the adjoint accurately. $\vect{R_U}$ and $\vect{J_U}$ are computed using automatic differentiation with CoDiPack \cite{codipack} and Trilinos’ Sacado package \cite{sacado} respectively.

To ensure that the measure of functional error we use for optimization is differentiable and to avoid the cancellations of positive and negative errors, we define the following measure for the goal-oriented error:
\begin{equation}
\label{F}
    \mathcal{F}(\vect{U},\vect{x}) = \frac{1}{2}\sum\limits_{i=1}^{n_{\text{dofs}}} \left(\psi_iR_i\right)^2 =\frac{1}{2} \vect{\eta^T\eta},
\end{equation}
where $\vect{\eta}$ is a vector of size $n_{\text{dofs}}$ and is computed as $\eta_i = \psi_iR_i$. Note that repeated indices do not imply summation in this work. While it is standard to use the cellwise dual-weighted residual, $|\vect{\psi_{\kappa}^TR_{\kappa}} |$, for fixed-fraction mesh adaptation, we observed empirically that the objective function in Eq. \ref{F} performed better for goal-oriented shock tracking. 

We note some distinguishing features of Eq. \ref{F} in comparison to the objective functions used in previous residual-based shock tracking approaches \cite{mdg1,zahr2}. In the MDG-ICE approach of Corrigan et al. \cite{mdg1}, the DG residual was formed without requiring a numerical flux at the face and the residual was augmented with interface conditions. This resulted in an overdetermined residual, which was solved in a least squares sense using unconstrained optimization. On the other hand, the shock tracking method of Zahr et al. \cite{zahr2} employed enriched residual to detect shocks using the idea that an enriched test space detects regions in the domain where the governing equation is violated, even though the PDE is orthogonal to the original test space. We note that in Eq. \ref{F}, the enriched DG residual is weighted with the adjoint solution. Adjoints indicate the sensitivity of a functional to residual perturbations in the domain and hence, indicate the importance of various regions in the domain for evaluating a functional. The objective function based on the adjoint-weighted residual, Eq. (\ref{F}), achieves goal-oriented shock tracking because $R_i$ is large in regions with discontinuities while the adjoint, $\psi_i$, is large in regions important for the functional. When the mesh faces are aligned with the shock, the solution within an element is smooth, resulting in a lower value of the objective function. Hence, the optimizer tries to align the mesh with the shock in regions deemed important for the functional to minimize the objective function based on the adjoint-weighted residual.

\subsection{First-order derivatives of the goal-oriented error indicator}
\label{first_order_derivatives}
Differentiating Eq. (\ref{F}) with respect to $\vect{x}$ and $\vect{U}$ gives
\begin{equation}
\label{F_derivatives}
    \begin{split}
        &\mathcal{F}_{\vect{x}} = \vect{\eta^T\eta_x}, \\
        &\mathcal{F}_{\vect{U}} = \vect{\eta^T\eta_U},
    \end{split}
\end{equation}
with $\vect{\eta_x}$ and $\vect{\eta_U}$ obtained using the chain rule
\begin{equation}
\label{eta_derivatives}
    \begin{split}
        &\vect{\eta_x} = \vect{\eta_\psi} \vect{\psi_x} + \vect{\eta_R}\vect{R_x} = \text{diag}(\vect{R})\vect{\psi_x} + \text{diag}(\vect{\psi})\vect{R_x}, \\
        &\vect{\eta_U} = \vect{\eta_\psi} \vect{\psi_U} + \vect{\eta_R}\vect{R_U} = \text{diag}(\vect{R}) \vect{\psi_U} + \text{diag}(\vect{\psi})\vect{R_U}.
    \end{split}
\end{equation}
In Eq. (\ref{eta_derivatives}), we have substituted $\vect{\eta_\psi} = \text{diag}(\vect{R})$ i.e. $[\vect{\eta_\psi}]_{ii} = R_i,\;i=1,..,n_{\text{dofs}}$. Similarly, $\vect{\eta_R} = \text{diag}(\vect{\psi})$. The derivatives of the adjoint can be obtained by differentiating the adjoint equation, Eq. (\ref{adjoint_eq}), to obtain
\begin{equation}
\label{adjoint_derivatives}
    \begin{split}
        &\vect{\psi_x} = \vect{-R_U^{-T}}\left(\vect{J_{Ux}} + \vect{\psi^TR_{Ux}}\right), \\
        &\vect{\psi_U} = \vect{-R_U^{-T}}\left(\vect{J_{UU}} + \vect{\psi^TR_{UU}}\right).
    \end{split}
\end{equation}
The second-order derivatives on the right-hand side of Eq. (\ref{adjoint_derivatives}) are computed using automatic differentiation.
\newdefinition{remark}{Remark}
\begin{remark}
\label{rmk1}
In order to evaluate $\vect{\psi^TR_{Ux}}$ and $\vect{\psi^TR_{UU}}$ without storing third-order tensors, we employ the following steps:
\begin{enumerate}
    \item Precompute and store the adjoint solution, $\vect{\psi}$, using Eq. (\ref{adjoint_eq}).
    \item Evaluate the scalar quantity $a = \vect{\psi^TR(U,x)}$ by declaring $\vect{U}$ and $\vect{x}$ as automatic differentiation variables.
    \item Automatically differentiate the scalar $a = \vect{\psi^TR(U,x)}$ twice with respect to $\vect{U}$ and $\vect{x}$.
\end{enumerate}
\end{remark}
This idea is drawn from a similar procedure employed in the full space aerodynamic shape optimization of Shi-Dong and Nadarajah \cite{doug} to evaluate $\vect{\lambda^Tr_{xx}}$, $\vect{\lambda^Tr_{xu}}$ and $\vect{\lambda^Tr_{uu}}$ for an arbitrary $\vect{\lambda}$ in the Hessian of the KKT system. Note that steps 2 and 3 can be performed cellwise, with their contributions added to the global second-order tensor. Finally, note that we never store the derivatives in Eqs. (\ref{eta_derivatives}) and (\ref{adjoint_derivatives}). Instead, we compute the products of $\vect{\psi_x}$, $\vect{\psi_U}$, $\vect{\eta_x}$ and $\vect{\eta_U}$ with vectors on the fly. Linear solves with $\vect{R_u^{-T}}$ in the products of $\vect{\psi_x}$ and $\vect{\psi_U}$ with vectors are evaluated at a tolerance of machine precision.

%%%%%%%%%%%%%%%%%%%%%%%%%%%%%%%%%%%
%%%%%%%%%%%%%%%%%%%%%%%%%%%%%%%%%%%%
\subsection{Second-order derivatives of the goal-oriented error indicator}
Second-order information about the goal-oriented error is essential to accelerate the convergence of the optimizer. Hence, we now derive the second-order derivatives $\mathcal{F}_{\vect{xx}}$, $\mathcal{F}_{\vect{Ux}}$ and $\mathcal{F}_{\vect{UU}}$. Differentiating Eq. (\ref{F_derivatives}) yields
\begin{equation}
\label{F_second_derivatives}
\begin{split}
       &    \mathcal{F}_{\vect{xx}} = \vect{\eta_x^T\eta_x} + \vect{\eta^T\eta_{xx}}, \\
       &   \mathcal{F}_{\vect{Ux}} = \vect{\eta_U^T\eta_x} + \vect{\eta^T\eta_{Ux}}, \\
       &   \mathcal{F}_{\vect{UU}} = \vect{\eta_U^T\eta_U} + \vect{\eta^T\eta_{UU}}.
\end{split}
\end{equation}
The first term on the right-hand side of Eq. (\ref{F_second_derivatives}) was discussed in the previous section. In order to form the second term on the right-hand side without storing third-order tensors, the procedure for forming $\vect{\eta^T\eta_{Ux}}$ is shown. The remaining terms can be derived similarly and their final forms are given at the end of this section.

We begin by noting that
\begin{equation}
\label{equality}
    \vect{\eta^T\eta_{Ux}} = \left(\vect{\Hat{\eta}^T\eta_U}\right)_{\vect{x}},
\end{equation}
where $\vect{\Hat{\eta}} = \vect{\eta}$, with the superscript ( $\Hat{}$ ) implying that $\vect{\Hat{\eta}}$ is being treated as a constant with respect to differentiation. Using the expression for $\vect{\eta_U}$ from Eq. (\ref{eta_derivatives}) gives
\begin{equation}
\label{eq_intermediate}
    \vect{\Hat{\eta}^T\eta_U} = \vect{\Hat{\eta}^T}\text{diag}(\vect{R})\vect{\psi_U} + \vect{\Hat{\eta}^T}\text{diag}(\vect{\psi}) \vect{R_U} = \vect{g^T\psi_U} + \vect{z^TR_U},
\end{equation}
with $\vect{g} =  \text{diag}(\vect{R})\vect{\Hat{\eta}}$, $g_i = \Hat{\eta}_{i}R_i$ and $\vect{z} = \text{diag}(\vect{\psi})\vect{\Hat{\eta}}$, $z_i = \Hat{\eta}_{i}\psi_i$. 

Differentiating Eq. (\ref{eq_intermediate}) with $\vect{x}$ and using the equality in Eq. (\ref{equality}) gives
\begin{equation}
\label{eta_uux_before}
     \vect{\eta^T\eta_{Ux}} = \vect{\psi_U^Tg_x} + \vect{g^T\psi_{Ux}} + \vect{R_U^Tz_x} + \vect{z^TR_{Ux}}.
\end{equation}
Since $\vect{\Hat{\eta}}$ is being treated as a constant, we have
\begin{equation}
\label{g_x_z_x}
    \begin{split}
        & \vect{g_x} = \vect{g_R}\vect{R_x} = \text{diag}(\vect{\Hat{\eta}})\vect{R_x}, \\
        & \vect{z_x} = \vect{z_\psi}\vect{\psi_x} = \text{diag}(\vect{\Hat{\eta}})\vect{\psi_x} .
    \end{split}
\end{equation}
Note that $g_i = \Hat{\eta}_{i}R_i = \psi_i R_i^2$ and, assuming that the residual $\vect{R}$ in the fine space computed from the interpolated solution will be small in magnitude as the optimizer converges, we assume that the term  $\vect{g^T\psi_{Ux}}$ in Eq. (\ref{eta_uux_before}) will be small in magnitude in comparison to the other terms. Thus, we neglect it and avoid computing the non-trivial term $\vect{\psi_{Ux}}$. As our test cases indicate, neglecting the term still retains the quadratic convergence of the optimizer. With these changes, substituting Eq. (\ref{g_x_z_x}) into Eq. (\ref{eta_uux_before}) and replacing $\vect{\Hat{\eta}}$ by $\vect{\eta}$ gives
\begin{equation}
    \label{eta_eta_ux}
      \vect{\eta^T\eta_{Ux}} \approx \vect{\psi_U^T}\text{diag}(\vect{\eta})\vect{R_x} + \vect{R_U^T}\text{diag}(\vect{\eta})\vect{\psi_x} + \vect{z^TR_{Ux}}.
\end{equation}
Using a similar procedure, we obtain the remaining terms of Eq. (\ref{F_second_derivatives}):
\begin{equation}
    \label{eta_eta_xx}
      \vect{\eta^T\eta_{xx}} \approx \vect{\psi_x^T}\text{diag}(\vect{\eta})\vect{R_x} + \vect{R_x^T}\text{diag}(\vect{\eta})\vect{\psi_x} + \vect{z^TR_{xx}},
\end{equation}
\begin{equation}
    \label{eta_eta_uu}
      \vect{\eta^T\eta_{UU}} \approx \vect{\psi_U^T}\text{diag}(\vect{\eta})\vect{R_U} + \vect{R_U^T}\text{diag}(\vect{\eta})\vect{\psi_U} + \vect{z^TR_{UU}},
\end{equation}
with 
\begin{equation}
    \label{z_eq}
    \vect{z} = \text{diag}(\vect{\psi})\vect{\eta}.
\end{equation}
The second derivatives of $\mathcal{F}$ can be computed by using Eqs. (\ref{eta_derivatives}), (\ref{eta_eta_ux}), (\ref{eta_eta_xx}), (\ref{eta_eta_uu}) and (\ref{z_eq}) in Eq. (\ref{F_second_derivatives}). It is not required to store the second derivatives of $\mathcal{F}$. Instead, we compute their products with vectors while using Krylov solvers to solve the resulting KKT system in optimization. Finally,  $\vect{z^TR_{xx}}$, $\vect{z^TR_{Ux}}$ and $\vect{z^TR_{UU}}$ are computed as described in Remark \ref{rmk1}.

\begin{comment}
Finally, we note that the linear solve with $\vect{R_U}$ required in the vector product with $\vect{\psi_x}$ and $\vect{\psi_U}$ might be expensive, especially if the number of vector products to be evaluated are large. In order to reduce the cost of optimization, we use the exact $\vect{R_U}$ while evaluating the gradient, eq. \ref{F_derivatives}, and use the preconditioner $\vect{\tilde{R}_U}$ in the vector products with the second-order derivatives, eq. \ref{F_second_derivatives}.  
\end{comment}

%% file: 5FullSpaceOptimizer.tex
\section{Full space optimizer for goal-oriented shock tracking}
\label{full_space}
We perform goal-oriented shock tracking by using an optimizer to minimize the goal-oriented error indicator, $\mathcal{F}(\vect{U},\vect{x})$. In order to ensure that the optimal mesh is valid, the penalty from Zahr and Persson \cite{zahr1} is used where a term, $f_{msh}$, is added to the objective function to penalize mesh distortion. The objective function we use for optimization is
\begin{equation}
\label{F_augmented}
    \mathcal{\tilde{F}}(\vect{U},\vect{x}) = \mathcal{F}(\vect{U},\vect{x}) + \alpha f_{msh}(\vect{x}),
\end{equation}
with $\vect{U} = \vect{I_h}\vect{u}$ as described in section \ref{goal_oriented_derivatives} and $\alpha \in \mathbb{R}^+$. $f_{msh}$ is taken to be \cite{zahr1}
\begin{equation}
\label{objfunc}
    f_{msh} = \frac{1}{N_k} \sum_{\kappa \in \mathcal{T}_h} \frac{1}{\int_\kappa \vect{dx}} \int_\kappa \left(\frac{\|\vect{J}\|_{F}^2}{(\text{det}(\vect{J}))^{2/d}}\right)^2 \vect{dx},
\end{equation}
where $d>1$ is the dimension of the domain and $\|\cdot\|_F$ refers to the Forbenius norm. $f_{msh}$ is differentiated with respect to $\vect{x}$ using automatic differentiation. 

The control variables for optimization, $\vect{s}$, are taken to be a subset of mesh nodes, $\vect{x}$. For example, if we would like to constrain the movement of nodes along one of the coordinates at a boundary, we do not include the fixed coordinate of the node in our control variable, $\vect{s}$. In this work, we use a linear mapping between $\vect{s}$ and $\vect{x}$ with $\vect{x(s)} = \vect{A}\vect{s}$, where $\vect{A}$ is a constant matrix.

Using the objective function, Eq. (\ref{objfunc}), the optimization problem can be stated as
\begin{equation} \label{optproblem}
\begin{array}{rrclcl}
\displaystyle \min_{\vect{u},\vect{s}} & \mathcal{\tilde{F}}(\vect{U(u)},\vect{x(s)})\\
\textrm{s.t.} & \vect{r(\vect{u},\vect{x(s)}}) = \vect{0}.
\end{array}
\end{equation}
The first and second-order derivatives with respect to design variables can be computed using chain rule:
\begin{equation}
\label{F_tilde_chain_rule}
\setlength{\jot}{10pt}
\begin{split}
    &\mathcal{\tilde{F}}_{\vect{u}} = \mathcal{\tilde{F}}_{\vect{U}} \frac{\partial \vect{U}}{\partial \vect{u}} = \mathcal{\tilde{F}}_{\vect{U}} \vect{I_h},\\
    &\mathcal{\tilde{F}}_{\vect{s}} = \mathcal{\tilde{F}}_{\vect{x}} \frac{\partial \vect{x}}{\partial \vect{s}} = \mathcal{\tilde{F}}_{\vect{x}} \vect{A}, \\
    &\mathcal{\tilde{F}}_{\vect{uu}} = \vect{I_h^T} \mathcal{\tilde{F}}_{\vect{UU}} \vect{I_h}, \\
    &\mathcal{\tilde{F}}_{\vect{us}} = \vect{I_h^T} \mathcal{\tilde{F}}_{\vect{Ux}} \vect{A}, \\
    &\mathcal{\tilde{F}}_{\vect{ss}} = \vect{A^T} \mathcal{\tilde{F}}_{\vect{Ux}} \vect{A}, \\
    &\vect{r_s} = \vect{r_x A}, \\
    &\vect{\lambda^T r_{us}} = \vect{\lambda^T r_{ux}A}, \\
    &\vect{\lambda^T r_{ss}} = \vect{A^T\lambda^T r_{xx}A},
\end{split}
\end{equation}
where $\vect{\lambda}$ is an arbitrary vector. Note that the interpolation operator, $\vect{I_h}$, depends only on the reference coordinates and is therefore constant with respect to the control variables.

In the classical reduced space optimization, the constraint $\vect{r(u,x(s))} = \vect{0}$ is eliminated by using the implicit function theorem to give $\vect{u} = \vect{u(x(s))}$. While the reduced space approach simplifies the optimization formulation, the cost of optimization scales rapidly with problem size. Moreover, the reduced-space approach requires the solution $\vect{u(x(s))}$ to exist for all intermediate iterates of $\vect{s}$ as the optimizer converges. This might not be possible in test cases involving shocks. It was observed in \cite{doug,hicken1} that the full space optimization, which solves for optimality and flow constraints simultaneously, is significantly more efficient and requires less computational resources if a suitable preconditioner for the residual Jacobian is adopted. Hence, in this work, we use full space optimization to solve the problem (\ref{optproblem}). Note that in the current work, we use second-order derivatives of the objective function and the residual in Eq. \ref{F_tilde_chain_rule}, in contrast to the Gauss-Newton Hessian approximation used in \cite{mdg1,zahr2}.

We form the Lagrangian function using a Lagrange multiplier, $\vect{\lambda}$:
\begin{equation}
    \mathcal{L}(\vect{u},\vect{s},\vect{\lambda}) = \mathcal{\tilde{F}}(\vect{U(u)}, \vect{x(s)}) + \vect{\lambda^T}\vect{r(u,x(s))}. 
\end{equation}
Since the optimal design and simulation variables minimize $\mathcal{L}$, we can obtain a local minimum by solving for those design variables that satisfy 
\begin{equation}
\label{lagrangian_gradient}
    \vect{\nabla \mathcal{L}} = \begin{bmatrix}
    \vect{\nabla_u \mathcal{L}}\\
    \vect{\nabla_s \mathcal{L}}\\
    \vect{\nabla_{\lambda} \mathcal{L}} 
    \end{bmatrix} = 
    \begin{bmatrix}
        \mathcal{\tilde{F}}_{\vect{u}} + \vect{\lambda^Tr_u} \\
        \mathcal{\tilde{F}}_{\vect{s}} + \vect{\lambda^Tr_s} \\
        \vect{r}
    \end{bmatrix} =
    \vect{0}.
\end{equation}

Using Newton's method to solve Eq. (\ref{lagrangian_gradient}) gives the following Karush-Kuhn-Tucker (KKT) system:
\begin{equation}
\label{KKT}
\begingroup
\renewcommand*{\arraystretch}{1.25}
    \begin{bmatrix}
\mathcal{L}_{\vect{uu}} & \mathcal{L}_{\vect{us}} & \vect{r_u^T}            \\
\mathcal{L}_{\vect{su}} & \mathcal{L}_{\vect{ss}} & \vect{r_s^T}            \\
\vect{r_u} &
\vect{r_s} &
\vect{0}
    \end{bmatrix}
\begin{bmatrix}
    \vect{p_u} \\
    \vect{p_s} \\
    \vect{p_\lambda}
\end{bmatrix}
    =
    \vect{-} \begin{bmatrix}
        \mathcal{\tilde{F}}_{\vect{u}} + \vect{\lambda^Tr_u} \\
        \mathcal{\tilde{F}}_{\vect{s}} + \vect{\lambda^Tr_s} \\
        \vect{r}
    \end{bmatrix},
\endgroup
\end{equation}
which is solved at each step of optimization to obtain the search direction. 

The full space KKT system, Eq. (\ref{KKT}), is often ill-conditioned and requires suitable preconditioners. Biros and Ghattas \cite{birosghattas1} derived a set of preconditioners based on the reduced space SQP, which was rigorously investigated in \cite{doug}. It was found in \cite{doug} that the $\vect{\tilde{P}_2}$ preconditioner of \cite{birosghattas1} was efficient for aerodynamic shape optimization. Hence, we employ the $\vect{\tilde{P}_2}$ preconditioner in this work, which can be written as \cite{birosghattas1}:
\begin{equation}
\label{P2A_preconditioner}
   \vect{\tilde{P}_2}   =  \begin{bmatrix}
            \mathbf{0} & \mathbf{0} & \mathbf{I} \\
            \mathbf{0} & \mathbf{I} & \mathbf{\mathbf{r_x}}^T\mathbf{\mathbf{\tilde{r}_u}}^{-T} \\
             \mathbf{I} & \mathbf{0} & \mathbf{0}
        \end{bmatrix}
               \begin{bmatrix}
            \mathbf{\tilde{r}_u} & \mathbf{r_x} & \mathbf{0} \\
            \mathbf{0} & \mathbf{B_z} & \mathbf{0} \\
             \mathbf{0} & \mathbf{0} & \mathbf{\tilde{r}_u}^{T}
        \end{bmatrix}, 
\end{equation}
where $\vect{\tilde{r}_u}$ is the preconditioner of the residual Jacobian and $\vect{B_{z}}$ is the approximation of the reduced space Hessian. Similar to \cite{doug}, we use Broyden-Fletcher-Goldfarb-Shanno (BFGS) quasi-Newton approximation of $\vect{B_z}$. $\vect{\tilde{P}_2}$ preconditioner approximates the KKT system matrix and is obtained by factoring the KKT system matrix on the left-hand-side of Eq. (\ref{KKT}), followed by replacing the instances of $\vect{r_u}$ and the reduced space Hessian by their approximations \cite{birosghattas1}. Preconditioning is achieved by applying $\vect{\tilde{P}_2^{-1}}$ on either sides of Eq. (\ref{KKT}). An expression for $\vect{\tilde{P}_2^{-1}}$ can be found in \cite[Appendix C]{doug}. In this work, Eq. (\ref{KKT}) is solved using FGMRES \cite{fgmres} with a relative tolerance of $10^{-6}$. 

While Newton's method, Eq. (\ref{KKT}), works well when the initial design variables are within the ball of convergence of the optimizer, it requires globalization strategies for convergence when started from an arbitrary initial guess. In the following section, we discuss the strategies employed to converge the full space optimizer in this work.

\subsection{Globalization strategies for full space optimization}

\subsubsection{Backtracking line search}
After solving Eq. (\ref{KKT}) for the search direction, the variables are updated with a step length $\gamma \in \mathbb{R}^+$, yielding $\vect{u}^{n+1} = \vect{u}^n + \gamma\vect{p_u}^n$, $\vect{s}^{n+1} = \vect{s}^n + \gamma\vect{p_s}^n$ and $\vect{\lambda}^{n+1} = \vect{\lambda}^n + \gamma\vect{p_\lambda}^n$. While the reduced space optimization chooses $\gamma$ by performing backtracking to ensure a sufficient decrease in the objective function, backtracking in full space optimization needs to reduce the objective function while ensuring that the flow constraint is not significantly violated. To that end, we choose the step length $\gamma$ to ensure sufficient decrease in the augmented Lagrangian \cite{doug}:
\begin{equation}
\label{lagrangian_augmented}
    \mathcal{L}^+(\vect{u},\vect{s},\vect{\lambda}) = \mathcal{L}(\vect{u},\vect{s},\vect{\lambda}) + \frac{\hat{\mu}}{2} \vect{r^Tr},
\end{equation}
with $\hat{\mu} = \max \{1,\; \mu/\| \mathcal{L}_{\vect{s}}\|\}$ and $\mu \in \mathbb{R}^+$. Note that the merit function heavily penalizes the violation of the flow constraint near optimality. The augmented Lagrangian is known to yield a step length of 1 as the optimizer converges \cite{birosghattas2, doug}, retaining the quadratic convergence of Newton's method.  

For the purpose of shock tracking, we also make sure that the updated mesh, $\vect{x}^{n+1} = \vect{x}(\vect{s}^{n+1})$, is valid. If the step length $\gamma$ violates the condition $\text{det}(\vect{J}^{n+1}) > 0$ in an element, $\gamma$ is reduced until the mesh is well-defined.

\subsubsection{Addition of regularization to the Hessian}
\label{section_regularization}
In order to control the magnitude of the updates and to improve conditioning of the KKT system, we use an approach of regularization on $\mathcal{L}_{\vect{ss}}$ similar to Zahr et al.\cite{zahr2}. A symmetric positive definite matrix $\varepsilon_s \vect{D}$ is added to $\mathcal{L}_{\vect{ss}}$, where $\vect{D}$ is either an identity or the stiffness matrix of the Poisson's equation and $\varepsilon_s \in \mathcal{R}^+$. Unless otherwise specified, $\vect{D}=\vect{I}$ in the test cases. Further, it was observed that in few of the test cases, the solution can become oscillatory as the mesh tracks the shock. Hence, to reduce oscillations in the solution, we also regularize $\mathcal{L}_{\vect{uu}}$ by adding $\varepsilon_u \vect{I}$, resulting in the following system:
\begin{equation}
\label{regularized_system}
\begingroup
\renewcommand*{\arraystretch}{1.25}
    \begin{bmatrix}
\mathcal{L}_{\vect{uu}} + \varepsilon_u \vect{I} & \mathcal{L}_{\vect{us}} & \vect{r_u^T}            \\
\mathcal{L}_{\vect{su}} & \mathcal{L}_{\vect{ss}} + \varepsilon_s \vect{D} & \vect{r_s^T}            \\
\vect{r_u} &
\vect{r_s} &
\vect{0}
    \end{bmatrix}
\begin{bmatrix}
    \vect{p_u} \\
    \vect{p_s} \\
    \vect{p_\lambda}
\end{bmatrix}
    =
    \vect{-} \begin{bmatrix}
        \mathcal{\tilde{F}}_{\vect{u}} + \vect{\lambda^Tr_u} \\
        \mathcal{\tilde{F}}_{\vect{s}} + \vect{\lambda^Tr_s} \\
        \vect{r}
    \end{bmatrix},
\endgroup
\end{equation}
where $\varepsilon_u$ and $\varepsilon_s$ are  both adapted according to $\|\vect{p_s}\|$ using an equation similar to \cite{zahr2}:
\begin{equation}
\label{varepsilon_update}
    \varepsilon_{u/s\;n+1} = 
    \begin{cases}
    \frac{\varepsilon_{u/s\;n}}{2}, & \|\vect{p_s}\| < 0.01 \\
    2\varepsilon_{u/s\;n}, & \|\vect{p_s}\| > 0.1 \\
    \varepsilon_{u/s\;n}, & \text{otherwise}
    \end{cases}.
\end{equation}
With this update, when the search direction $\| \vect{p_s}\|$ is large, the regularization on the control ($\varepsilon_s$) is increased. Similarly, the idea of updating $\varepsilon_u$ in accordance with $\| \vect{p_s}\|$ is that, when the mesh moves towards the shock, it is desirable to have large regularization on the solution update to reduce Gibbs oscillations and when the mesh is almost stationary, $\varepsilon_{u/s} \to 0$ and the optimizer utilizes the efficient search direction from Newton's update.

%% file: 6Algorithm.tex
\subsection{Algorithm for full space goal-oriented implicit shock tracking}
\label{algorithm}
We now summarize the approach of full space goal-oriented shock tracking, as described in sections \ref{DG} to \ref{full_space}, in algorithm \ref{algo_eq}.

\begin{algorithm}[H]
\caption{Full space optimization for adjoint-based goal-oriented shock tracking}
\label{algo_eq}
\begin{algorithmic}[1]

\State \textbf{Initialize:} control variable $\vect{s_0}$, initial mesh $\vect{x} = \vect{x(s_0)}$, $\varepsilon_s = \varepsilon_{s0}$, $\varepsilon_u = \varepsilon_{u0}$ 
\State Choose a functional, $J$, of interest.
\State If the test case permits, evaluate $\vect{u}$ such that $\vect{r(u,x)} = \vect{0}$.
\State Evaluate $\vect{r_u}$, $\vect{r_x}$.
\State Compute interpolated solution $\vect{U} = \vect{I_h u}$.
\State Evaluate $\vect{R(U,x)}$, $\vect{R_U}$, $\vect{R_x}$.
\State Evaluate $\vect{J_U}$, $\vect{J_{Ux}}$ and $\vect{J_{UU}}$
\State Evaluate the adjoint, $\vect{\psi}$, from the adjoint equation \Comment{Eq. (\ref{adjoint_eq})}.
\State Evaluate $\vect{\psi^T R_{UU}}$ and $\vect{\psi^T R_{Ux}}$  \Comment{see Remark \ref{rmk1}}

\State Evaluate $\mathcal{\tilde{F}}_{\vect{u}}$ and $\mathcal{\tilde{F}}_{\vect{s}}$ \Comment {Eqs. (\ref{F_tilde_chain_rule}) and (\ref{F_derivatives})}

\State \textbf{Initialize:} $\vect{\lambda} = \vect{r_u^{-T}}\mathcal{\tilde{F}}_{\vect{u}}$

\State Evaluate $\vect{\nabla \mathcal{L}}$ 

\While{$\| \vect{\nabla \mathcal{L}} \| > \text{tol}$}
\State Evaluate $\vect{\lambda^Tr_{uu}}$, $\vect{\lambda^Tr_{ux}}$ and $\vect{\lambda^Tr_{xx}}$
\State Evaluate $\vect{z} = \text{diag}(\vect{\psi}) \vect{\eta}$
\State Evaluate $\vect{z^TR_{UU}}$, $\vect{z^TR_{Ux}}$ and $\vect{z^TR_{xx}}$ 
\State Solve regularized KKT system using Krylov solvers. \Comment{Eq. (\ref{regularized_system})} 
\State Update $\varepsilon_s$ and $\varepsilon_u$  \Comment{Eq. (\ref{varepsilon_update})}

\State Perform backtracking using the augmented Lagrangian \Comment{Eq. (\ref{lagrangian_augmented})}

\State Update variables: $\vect{u} \gets \vect{u} + \gamma \vect{p_u}$, $\vect{s} \gets \vect{s} + \gamma \vect{p_s}$, $\vect{\lambda} \gets \vect{\lambda} + \gamma \vect{p_{\lambda}}$
\State Update mesh: $\vect{x} = \vect{x(s)}$
\State Perform steps 4 to 10
\State Evaluate $\vect{\nabla \mathcal{L}}$
\EndWhile
\end{algorithmic}
\end{algorithm}
A tolerance of $\text{tol} = 10^{-10}$ is used to determine convergence of the algorithm. Full space optimization is performed in a parallel optimization framework using Trilinos' Rapid Optimization Library (ROL) \cite{rol}.  

Despite a few similarities, we note some important differences between the currently developed adjoint-based approach and the full space shock tracking approach of Zahr et al. \cite{zahr2}. In addition to the fact that we use the adjoint-weighted enriched residual as the objective function instead of just the enriched residual, the SQP solver employed in \cite{zahr2} differs as well. \cite{zahr2} uses a Gauss-Newton approximation of the Hessian, which neglects the second-order derivatives of the objective function and the residual. Further, preconditioners for the SQP method in \cite{zahr2} are being developed \cite{zahr_preconditioner}. In the current work, we use a standard Newton's method by including the second-order information to improve convergence at an additional cost. Using the standard SQP method allows us to enhance the scalability of the algorithm by using the full space preconditioners developed in \cite{birosghattas1}, which have been extensively tested in the past for aerodynamic shape optimization \cite{doug}.

In the sections that follow, test cases using algorithm \ref{algo_eq} are presented to demonstrate functional-dependent shock tracking.

%% file: 7NumericalTestCases.tex
\section{Numerical test cases}
\label{test_cases}
In this section, we present a series of test cases to demonstrate shock tracking using a functional of interest. PDEs of advection and inviscid Euler equations are considered, with solutions containing straight and curved shocks. We also specify the numerical flux, the regularization, the weight on mesh distortion and the penalization used for backtracking.

\subsection{Linear advection with a straight shock}
Consider the steady-state linear advection equation in $u$ on $(x,y) \in \Omega = [-1,1] \times [0,1]$:
\begin{equation}
\label{eq_advection_straight_shock}
\begin{split}
        &\vect{\nabla} \cdot \left( \vect{\zeta} u \right) = 0 \quad \text{in} \quad \Omega, \\
        & u = H(x) = 
        \begin{cases}
            1, & x \geq 0 \\
            0, & x<0 
        \end{cases}
        \quad \text{on} \quad \Gamma_{-} = \{(x,y) \in \Gamma \mid \vect{\zeta}\cdot \vect{n} \leq 0\},
\end{split}
\end{equation}
where $\vect{\zeta} = [-0.3\;\;1]^T$ is the velocity vector and $\Gamma_{-}$ is the inflow boundary. The functional of interest is taken to be
\begin{equation}
\label{functional_straight_shock}
    J(u) = \int_{\Gamma_{+}} (1 + x)u \;d\Gamma,
\end{equation}
where $\Gamma_{+}$ is the outflow boundary. Since the domain is rectangular in the current case, it can be determined from the velocity vector $\vect{\zeta}$ that $\Gamma_{+}$ consists of the left and the top boundaries while $\Gamma_{-}$ consists of the bottom and the right boundaries.

The boundary condition in Eq. (\ref{eq_advection_straight_shock}) is discontinuous at $(0,0)$, which creates a shock in $u$ starting at $(0,0)$ and propagating along the velocity $\vect{\zeta}$ throughout the domain until the outflow boundary. Therefore, to evaluate the functional $J$ at outflow in Eq. \ref{functional_straight_shock} accurately, it is expected that the shock in the entire domain needs to be well resolved.

\begin{figure}[H]
    \centering
    \includegraphics[scale=0.4]{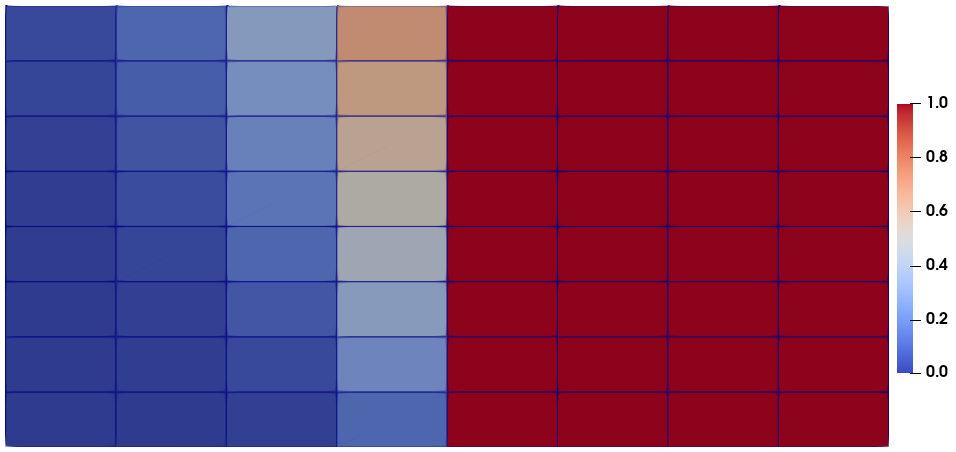}
    \caption{Initial mesh and solution.}
   \label{initial_mesh_straight_shock}
\end{figure}

We use $p=0$ order of approximation for the solution and $q=1$ basis functions for the grid. Upwind flux is used at elements' faces in the DG discretization of Eq. (\ref{eq_advection_straight_shock}). Fig. \ref{initial_mesh_straight_shock} shows the initial mesh and the computed $p=0$ solution.

As seen in Fig. \ref{initial_mesh_straight_shock}, the shock is not well captured by the structured grid. The optimizer is run as described in algorithm \ref{algo_eq}, with the functional in Eq. (\ref{functional_straight_shock}) and the DG discretization of Eq. (\ref{eq_advection_straight_shock}). We use $\mu = 0.01$ in the augmented Lagrangian for backtracking, remove the weight on mesh distortion i.e. $\alpha = 0$, and set  $\varepsilon_{s0} = 0.01$ and $\varepsilon_{u0} = 0$ for regularization. The control variables, $\vect{s}$, are the set of all nodes in the domain, with the nodes at the boundary constrained to slide along the boundary. Similar to the test case in \cite{zahr2}, the node at $(0,0)$ is fixed to ensure that the Heaviside boundary condition is evaluated and integrated accurately.

\begin{figure}[H]
    \centering
    \includegraphics[scale=0.4]{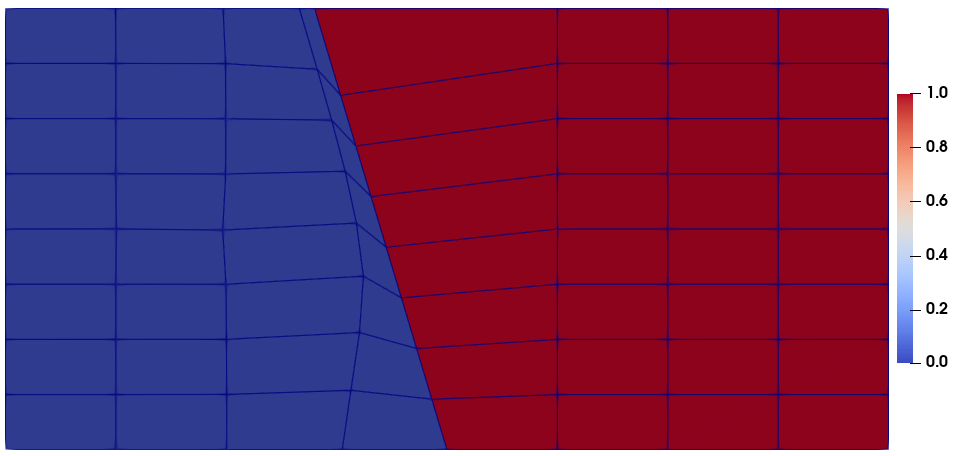}
    \caption{Converged mesh and solution from goal-oriented shock tracking.}
    \label{optimal_mesh_straight_shock}
\end{figure}

\begin{figure}[H]
    \centering
    \begin{subfigure}{0.55\textwidth}
    \includegraphics[width=1\linewidth]{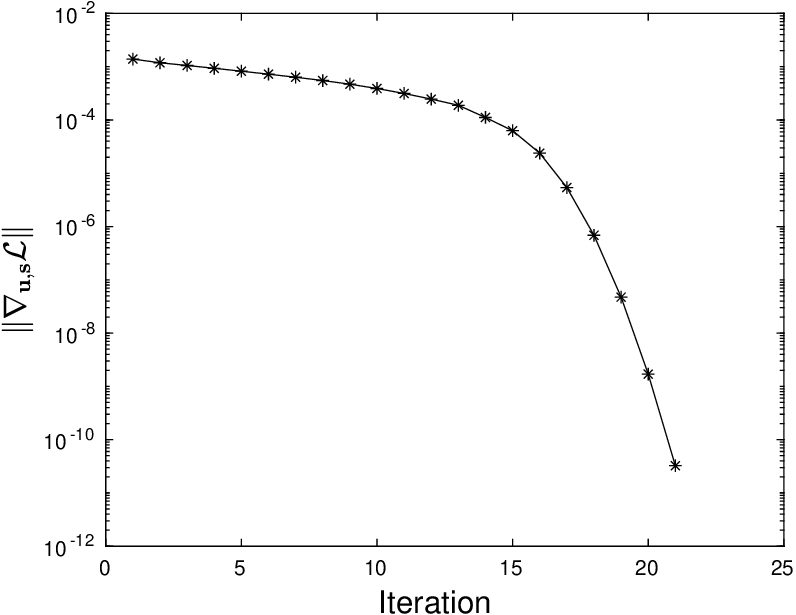}
    \caption{}
    \label{straight_shock_gradient_a}
    \end{subfigure}
    
    \begin{subfigure}{0.55\textwidth}
    \includegraphics[width=1\linewidth]{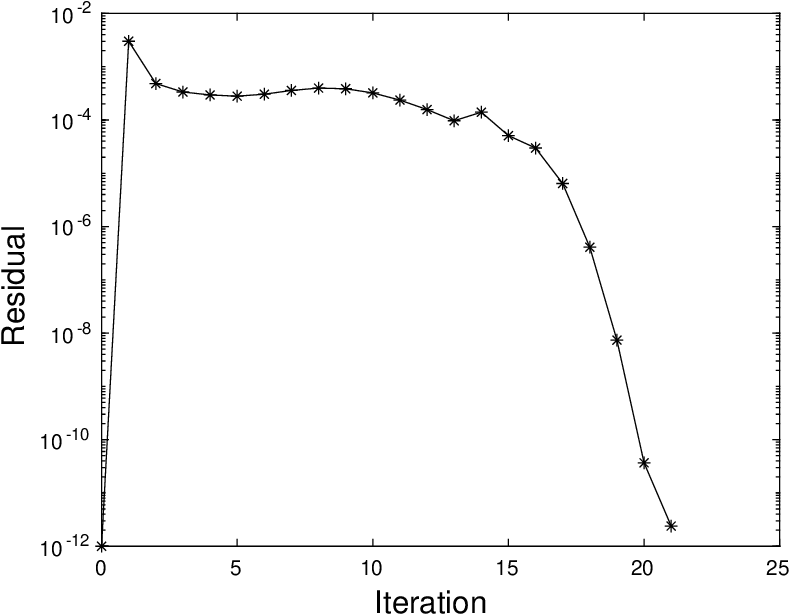}
    \caption{}
    \label{straight_shock_gradient_b}
    \end{subfigure}
    \caption{Convergence of Lagrangian gradient, $\vect{\nabla \mathcal{L}}$, for goal-oriented shock tracking. (a) convergence of $\|\vect{\nabla_{u,s} \mathcal{L}} \|$. (b) convergence of $\|\vect{\nabla_{\lambda}\mathcal{L}}\| =\|\vect{r(u,x)}\|$. }
    \label{lagrangian_gradient_straight_shock}
\end{figure}

\begin{figure}[H]
    \centering
    \includegraphics[scale=0.6]{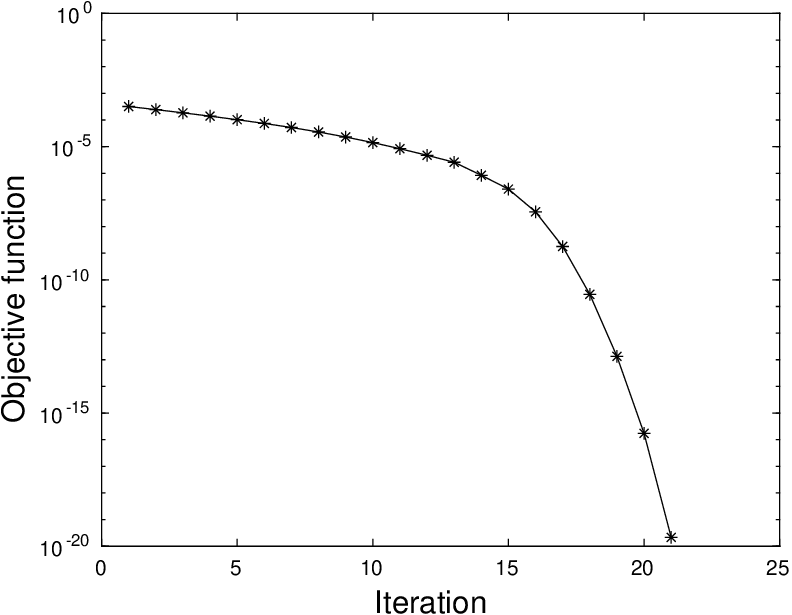}

\caption{Evolution of the objective function $\mathcal{\tilde{F}}$.}
\label{objective_function_straight_shock}
\end{figure}

The full space optimizer converges to the mesh and the solution as shown in Fig. \ref{optimal_mesh_straight_shock}. As expected, the optimal mesh which minimizes the functional error is aligned with the shock in the entire domain. Fig. \ref{straight_shock_gradient_a} shows the convergence of the gradient, $\vect{\nabla_{u,s}\mathcal{L}} = [\vect{\nabla_u \mathcal{L}} \quad \vect{\nabla_s \mathcal{L}}]^T$ and Fig. \ref{straight_shock_gradient_b} shows the convergence of the residual, $\vect{\nabla_{\lambda} \mathcal{L}} = \vect{r(u,x(s))}$. Quadratic convergence of Newton's method is observed. Fig. \ref{objective_function_straight_shock} demonstrates the evolution of the objective function. The objective function and hence, the error in the functional are close to machine precision on the optimal mesh.   

\begin{figure}[H]
    \centering \includegraphics[scale=0.6]{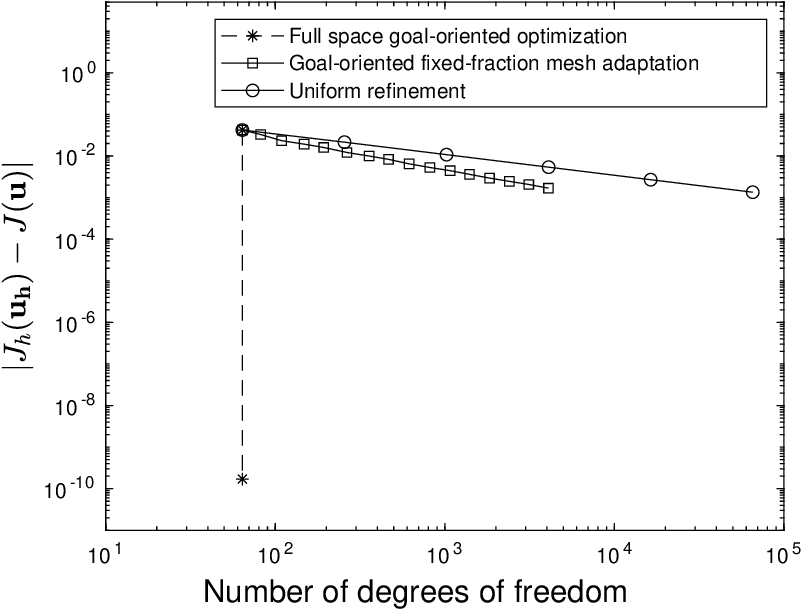}
\caption{Convergence of the functional error vs degrees of freedom for various mesh adaptation strategies.}
\label{functional_error_mesh_adaptations}
\end{figure}

\begin{figure}[H]
    \centering
    \includegraphics[scale=0.4]{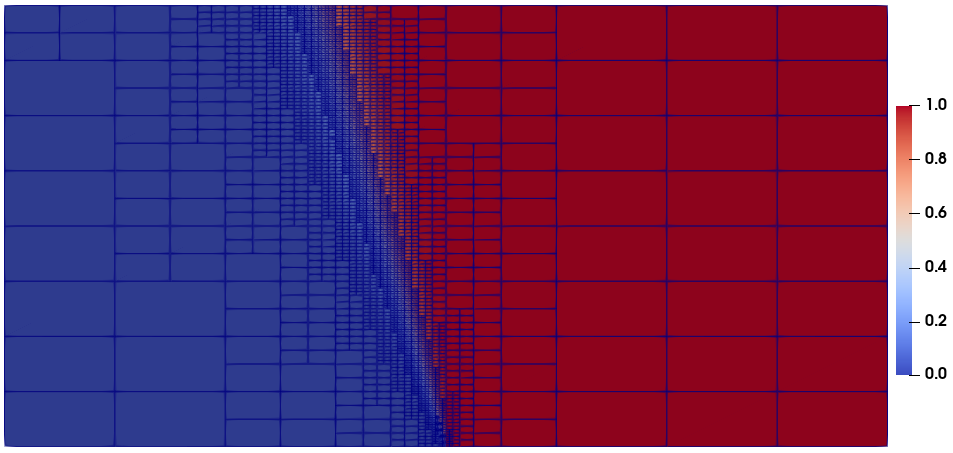}
    \caption{Mesh obtained from fixed-fraction goal-oriented $h$-adaptation.}
    \label{fixedfraction_straight_shock}
\end{figure}

The exact solution of Eq. \ref{eq_advection_straight_shock} is $u = H(x + 0.3y)$. Since the goal of performing optimization was to reduce error in the functional, we compare the error in evaluating the functional numerically, $| J_h(u_h) - J(u)|$, using the exact analytical value of the functional, $J(u) = 1.755$. Fig. \ref{functional_error_mesh_adaptations} compares the reduction in error in the functional using goal-oriented optimization, fixed-fraction based on refining $10\%$ of cells with the highest cellwise dual-weighted residual, and uniform $h$-refinement strategies. It can be observed that optimization achieves a significantly lower functional error when compared to other adaptation strategies while maintaining constant degrees of freedom. Fig. \ref{fixedfraction_straight_shock} shows the final adapted mesh obtained from goal-oriented fixed-fraction mesh adaptation after 15 adaptation cycles. 

\begin{figure}[H]
    \centering
    \includegraphics[scale=0.4]{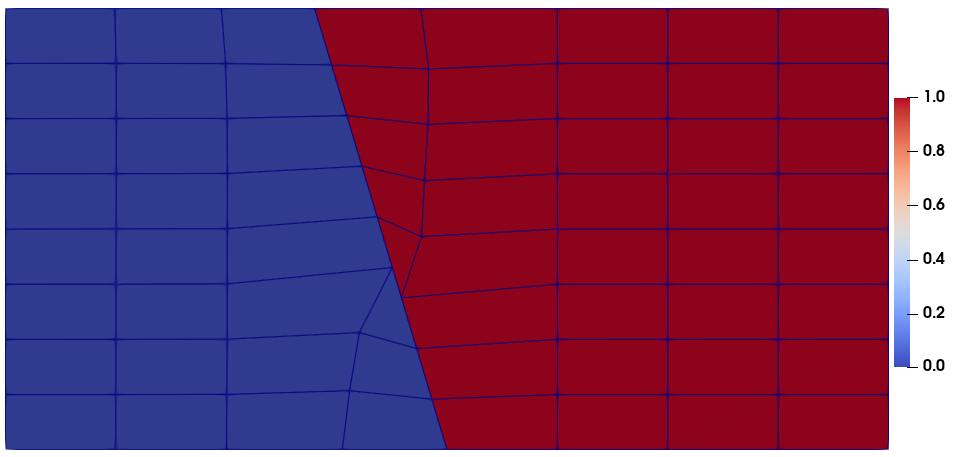}
    \caption{Converged mesh and solution from residual-based implicit shock tracking.}
    \label{optimal_mesh_straight_shock_zahr}
\end{figure}

\begin{figure}[H]
    \centering
    \begin{subfigure}{0.55\textwidth}
    \includegraphics[width=1\linewidth]{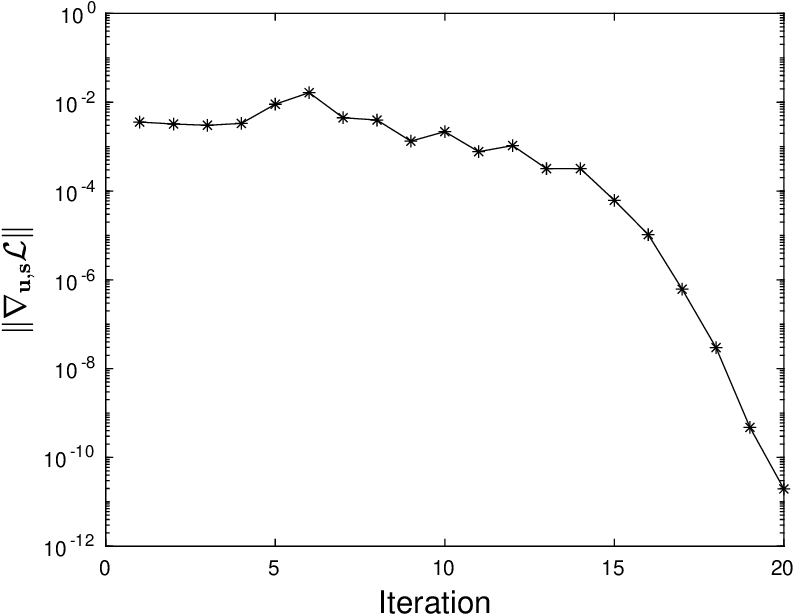}
    \caption{}
    \end{subfigure}
    % Add convergence of residual for zahr's approach here.
    \begin{subfigure}{0.55\textwidth}
    \includegraphics[width=1\linewidth]{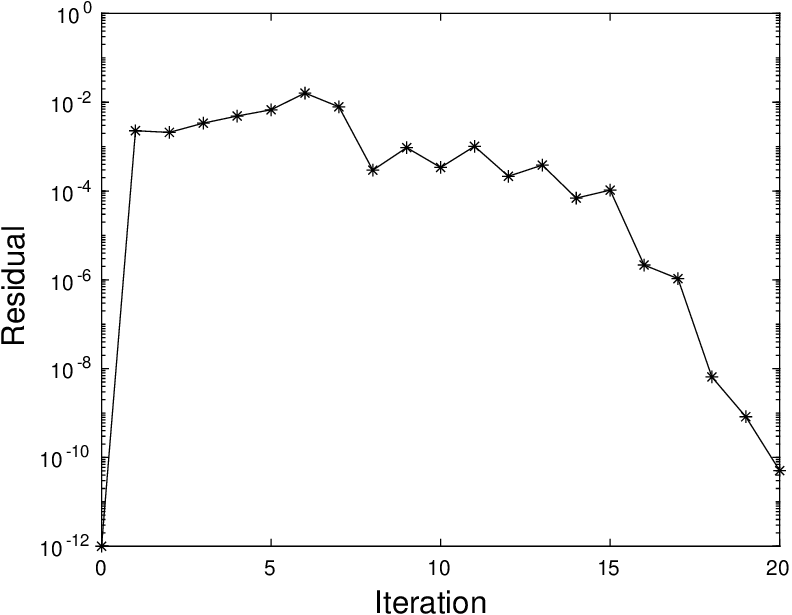}
    \caption{}
    \end{subfigure}
    \caption{Convergence of Lagrangian gradient, $\vect{\nabla \mathcal{L}}$, for residual-based shock tracking. (a) convergence of $\|\vect{\nabla_{u,s} \mathcal{L}} \|$. (b) convergence of $\|\vect{\nabla_{\lambda}\mathcal{L}}\| =\|\vect{r(u,x)}\|$. }
    \label{zahr_gradient_straight_shock}
\end{figure}

We also run our test case with the implicit residual-based shock tracking approach proposed in Zahr et al. \cite{zahr2}. Note that Zahr et al. \cite{zahr2} used an enriched test space with the original trial space for representing the solution, while we use an enriched test space and interpolate the solution to the enriched trial space to form the enriched residual $\vect{R}$.  Using the objective function $\vect{R^T}\vect{R}$ with the Newton SQP framework in section \ref{full_space}, we run the optimizer and obtain the optimal mesh as shown in Fig. \ref{optimal_mesh_straight_shock_zahr}. Fig. \ref{zahr_gradient_straight_shock} shows the convergence of the Lagrangian gradient and the residual for residual-based shock tracking.

%%%%%%%%%
%%%%%%%%%
Note that the goal-oriented shock tracking proposed in this work and the residual-based shock tracking both track the entire shock, which is expected when accurate evaluation of the functional requires the shock in the entire domain to be well-captured. The next test case demonstrates that adjoint-based shock tracking can track only a part of the shock if computing accurate functional requires only a subset of the solution field to be evaluated accurately.

%%%%%%%%%%%%%%%%%%%%%%%%%%%%%%%%%%%%%%%
%%%%%%%%%%%%%%%%%%%%%%%%%%%%%%%%%%%%%%%%
%%%%%%%%%%%%%%%%%%%%%%%%%%%%%%%%%%%%%%%%
\subsection{Vector advection equation with a curved shock}
\label{test_case_2}
Consider the set of advection equations on the domain $\Omega = [0\quad 1]^2$ in two states, $u_0$ and $u_1$:
\begin{equation}
\label{u0_eq_volume}
    \vect{\nabla} \cdot \left(\vect{\beta}u_0 \right)= 0 \quad \text{in} \quad \Omega,
\end{equation}
\begin{equation}
\label{u1_eq_volume}
    \vect{\nabla} \cdot \left(\vect{\mathcal{V}}u_1\right) = u_0^2 \quad \text{in} \quad \Omega, 
\end{equation}
with the boundary conditions:
\begin{equation}
\label{u0_eq_boundary}
    u_0 = \begin{cases}
        1 & x \leq 0.375 \\
        0 & x > 0.375
    \end{cases} \quad \text{on} \quad \Gamma_0=\{(x,y) \in \Gamma \mid \vect{\beta}\cdot \vect{n} \leq 0\},
\end{equation}

\begin{equation}
\label{u1_eq_boundary}
    u_1 = 0.3 \quad \text{on} \quad \Gamma_1=\{(x,y) \in \Gamma \mid \vect{\mathcal{V}}\cdot \vect{n} \leq 0\},
\end{equation}
and the velocity vectors:
\begin{equation}
\label{eq_velocity_vectors}
    \vect{\beta} = \begin{bmatrix}
        0.8y - 0.4 \\
        -1
    \end{bmatrix}, \quad \quad \vect{\mathcal{V}} = \begin{bmatrix}
        1 \\
        -1
    \end{bmatrix}.
\end{equation}
We are interested in the accurate evaluation of the following functional:
\begin{equation}
\label{eq_functional_case2}
    J(\vect{u}) = \int_{\Gamma_{\text{right}}} u_1^2 d\Gamma.
\end{equation}

\begin{figure}[H]
    \centering
    \includegraphics[scale=0.4]{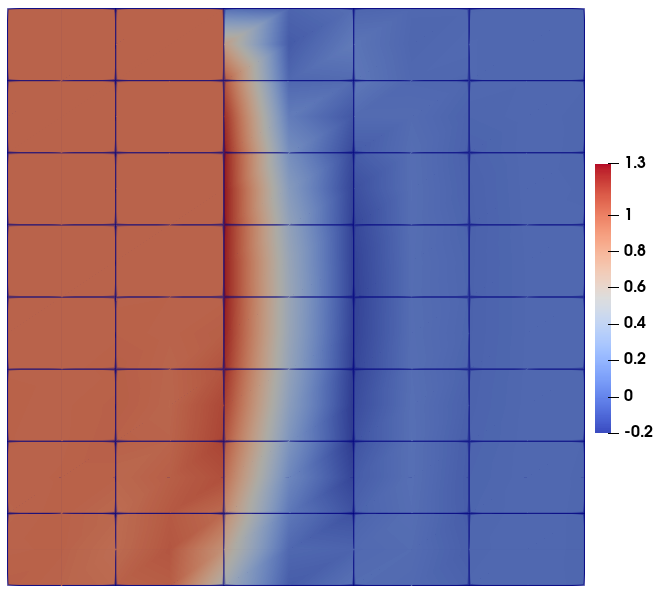}
    \caption{Initial mesh and solution $u_0$.}
    \label{initial_mesh_cshock}
\end{figure}

Before proceeding with the optimization, it is informative to understand the characteristics of the test case. Note that the state $u_0$ is completely specified by Eq. \ref{u0_eq_volume} and the boundary condition Eq. \ref{u0_eq_boundary}. Fig. \ref{initial_mesh_cshock} shows $u_0$ as computed on the initial mesh. The exact $u_0$ has a shock along $x = -0.4y^2 + 0.4y + 0.375$ and can be determined to be:
\begin{equation}
\begin{split}
    & u_0 = 1 \quad \text{in} \quad \Omega_1 = \{(x,y) \in \Omega \mid x \leq -0.4y^2 + 0.4y + 0.375 \} \\
     & u_0 = 0 \quad \text{in} \quad \Omega_2 = \{(x,y) \in \Omega \mid x > -0.4y^2 + 0.4y + 0.375 \}
\end{split} 
\end{equation}

\begin{figure}[H]
    \centering
    \includegraphics[scale=0.25]{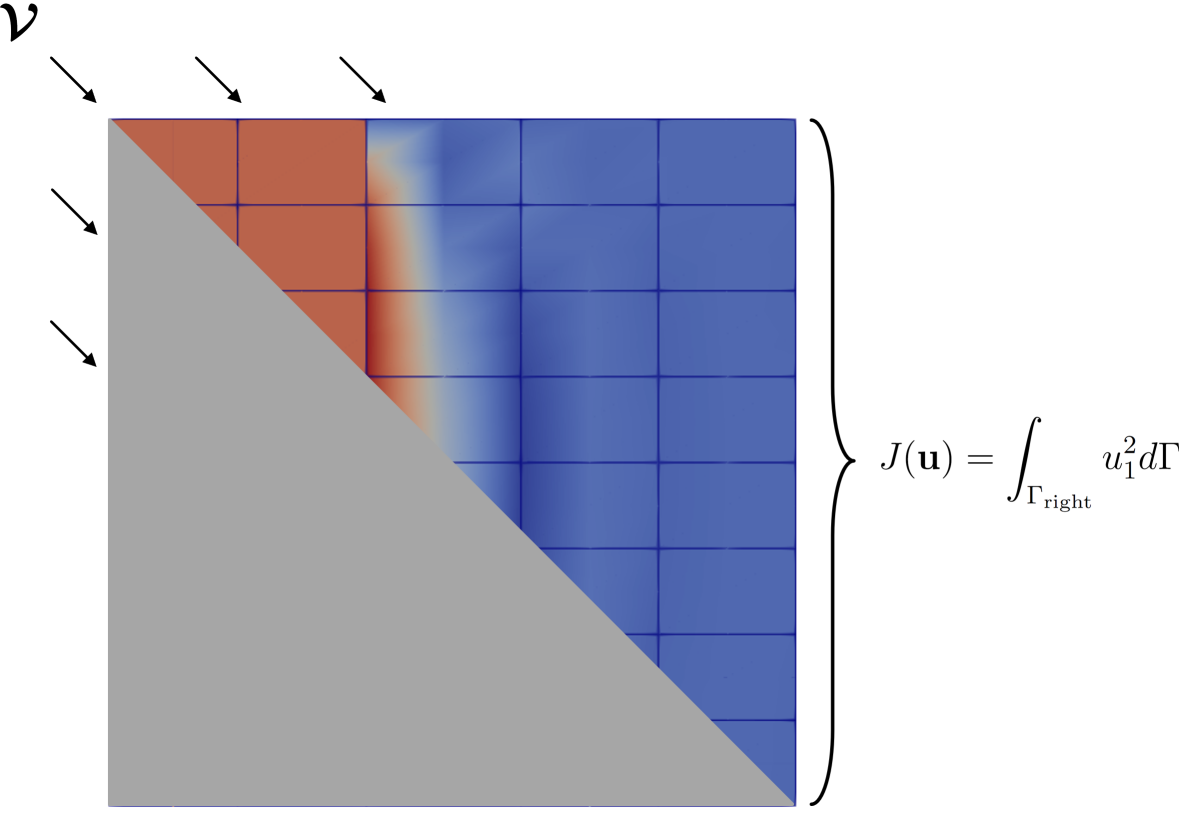}
    \caption{Region important for evaluating the functional. The arrows on top left indicate the direction of flow of state $u_1$.}
    \label{fig_directional_flow_u1}
\end{figure}

State $u_1$ depends on state $u_0$ through Eq. \ref{u1_eq_volume}. If $u_0$ is computed accurately, the equation for $u_1$ in Eq. \ref{u1_eq_volume} can be written as:
\begin{equation}
    \vect{\nabla} \cdot \left(\vect{\mathcal{V}} u_1\right) = 1 \quad \text{in} \quad \Omega_1 
\end{equation}
\begin{equation}
    \vect{\nabla} \cdot \left(\vect{\mathcal{V}} u_1\right) = 0 \quad \text{in} \quad \Omega_2 
\end{equation}
It can be seen that $u_1$ propagates along $\vect{\mathcal{V}}$ and is constant along $\vect{\mathcal{V}}$ in $\Omega_2$. Hence, the value attained by $u_1$ along $x = -0.4y^2 + 0.4y + 0.375$ propagates to the boundary $\Gamma_{\text{right}}$ along $\vect{\mathcal{V}}$, which is then used to evaluate the functional $J(\vect{u})$ in Eq. \ref{eq_functional_case2}. With $\vect{\mathcal{V}}$ in Eq. \ref{eq_velocity_vectors}, we note that only $u_1$ above the line $x+y=1$ is important to evaluate $J(\vect{u})$, as highlighted in Fig. \ref{fig_directional_flow_u1}. Finally, note that when $u_0$ is not evaluated accurately, the error in $u_0$ causes an error in $u_1$ and subsequently results in an error in evaluating $J(\vect{u})$.

To make the observation concrete and gain further insight into the adjoint-weighted residual objective function, we investigate the continuous adjoint equation related to the original PDE in Eqs. \ref{u0_eq_volume} - \ref{u1_eq_boundary} and functional in Eq. \ref{eq_functional_case2}. The derivation of the continuous adjoint equation can be found in~\ref{continuous_adjoint_appendix}, the final form of which is:
\begin{equation}
\label{eq_psi0_vol}
    \vect{\nabla}\psi_0 \cdot \vect{\beta} + 2u_0\psi_1 = 0 \quad \text{in} \quad \Omega,
\end{equation}
\begin{equation}
\label{eq_psi1_vol}
    \vect{\nabla}\psi_1 \cdot \vect{\mathcal{V}} = 0 \quad \text{in} \quad \Omega,
\end{equation}
with the adjoint boundary conditions:
\begin{equation}
\label{eq_psi0_boundary}
    \left(\vect{\beta}\cdot \vect{n}\right)\psi_0 = 0 \quad \text{on} \quad \Gamma \setminus \Gamma_0,
\end{equation}

\begin{equation}
\label{eq_psi1_boundary}
\left(\vect{\mathcal{V}}\cdot \vect{n}\right)\psi_1 = \begin{cases}
    2u_1, & \Gamma_{\text{right}} \\
    0, & \Gamma_{\text{bottom}}
\end{cases} .
\end{equation}
From Eq. \ref{eq_psi1_vol} for $\psi_1$, it can be noted that the boundary condition, Eq. \ref{eq_psi1_boundary}, is propagated parallel to $\vect{\mathcal{V}}$. Therefore, $\psi_1 = 0$ in the region below the line $x+y=1$ and the source term $2u_0\psi_1$ in Eq. \ref{eq_psi0_vol} is zero below the line $x+y=1$. As a result, from the boundary condition for $\psi_0$ in Eq. \ref{eq_psi0_boundary}, we can determine that $\psi_0 = 0$ below the line $x+y=1$. Hence, the adjoint solution, which is used to indicate regions important for the functional, is expected to indicate the region above $x+y=1$ as important for evaluating the functional. Therefore, we would expect the optimizer to track the shock in $u_0$ in this region. 

However, before proceeding with the optimization, note that $\psi_1$ has a shock due to the discontinuity in the boundary condition, Eq. \ref{eq_psi1_boundary}. Since we are differentiating the adjoint in our approach, it is of interest to use smooth adjoint solutions to enhance convergence of the optimizer. Hence, we modify the functional and use
\begin{equation}
\label{modified_functional_testcase2}
    \tilde{J}(\vect{u}) = \int_{\Gamma_{\text{right}}} g(y) u_1^2 d\Gamma, 
\end{equation}
leading to the adjoint boundary condition
\begin{equation}
\label{eq_psi1_smooth_boundary}
\left(\vect{\mathcal{V}}\cdot \vect{n}\right)\psi_1 = \begin{cases}
    2g(y)u_1, & \Gamma_{\text{right}} \\
    0, & \Gamma_{\text{bottom}}
\end{cases} .
\end{equation}

\begin{figure}[H]
    \centering
    \begin{subfigure}{0.55\textwidth}
    \includegraphics[width=1\linewidth]{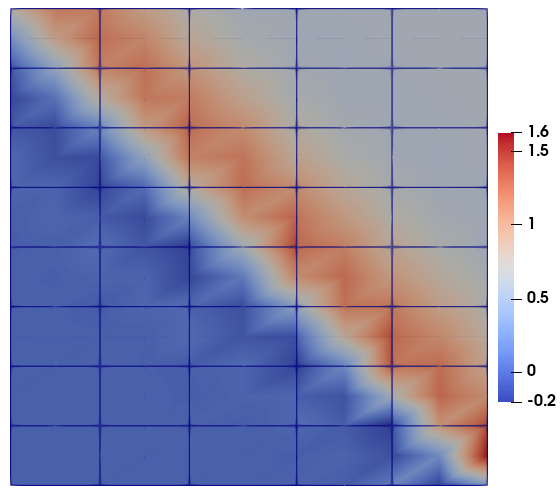}
    \caption{}
    \label{shocked_adjoint}
    \end{subfigure}
    
    \begin{subfigure}{0.55\textwidth}
    \includegraphics[width=1\linewidth]{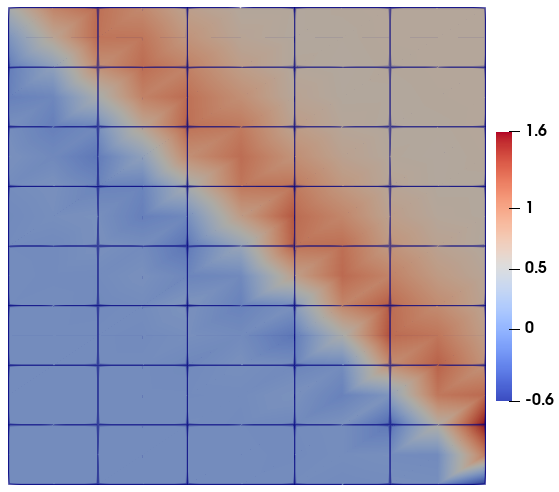}
    \caption{}
    \label{smoothened_adjoint}
    \end{subfigure}
    \caption{Shock in the adjoint solution, $\psi_1$, and making the adjoint smooth. (a) Adjoint solution using $J(\pmb{u}) = \int_{\Gamma_{\text{right}}}u_1^2 d\Gamma$. (b) Adjoint solution using $\tilde{J}(\pmb{u}) = \int_{\Gamma_{\text{right}}}g(y)u_1^2 d\Gamma$.}
    \label{shock_smooth_adjoint}
\end{figure}

In order to have a continuous $\psi_1$, $g(y)$ is chosen to be a continuous function in $y$ such that $g(y) \to 0$ at $\Gamma_{\text{right}} \cap \Gamma_{\text{bottom}}$ and $g(y) \to 1$ away from $\Gamma_{\text{bottom}}$. In this work, we use the logistic function
\begin{equation}
    g(y) = \frac{1}{1+e^{-230(y-0.05)}}.
\end{equation}
Fig. \ref{shock_smooth_adjoint} shows the shock in the adjoint solution with the original functional and the smooth adjoint solution with the modified functional. The discrete adjoint equation, Eq. \ref{adjoint_eq}, was used to obtain the coefficients of the adjoint solution $\vect{\psi}$ and, since we are using an adjoint consistent DG discretization, the discrete adjoint equation is a consistent discretization of the continuous adjoint problem, Eqs. \ref{eq_psi0_vol} - \ref{eq_psi1_boundary}, with $\vect{\psi} \in [\mathcal{P}^{p+1}(\hat{\kappa})]^2$. Finally, we observed that the optimization problem did not converge when using a non-smooth adjoint solution and it was critical to modify the functional before running the optimizer.

\begin{figure}[H]
    \centering
    \includegraphics[scale=0.4]{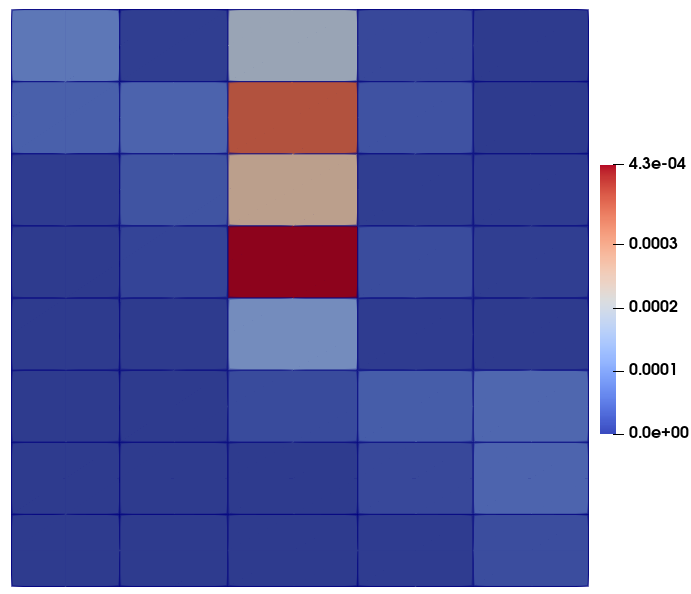}
    \caption{Cellwise dual-weighted residual error indicator on the initial mesh.}
    \label{dwr_cshock}
\end{figure}

Fig. \ref{dwr_cshock} shows the distribution of cellwise dual-weighted residual, $|\vect{\psi_{\kappa}^TR_{\kappa}}|$, highlighting that most of the error in the functional is due to the shock in $u_0$ in the top-half of the domain.  

\begin{figure}[H]
    \centering
    \includegraphics[scale=0.4]{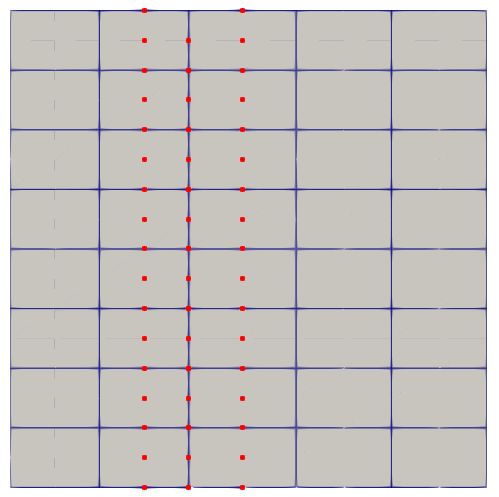}
    \caption{Control nodes used in optimization for grid of degree 2. Grid of degree 1 uses only the highlighted nodes at the vertices.}
    \label{control_nodes_cshock}
\end{figure}

We run the full space optimization using a solution of order $p=2$ and grids of orders 1 and 2. Fig. \ref{initial_mesh_cshock} shows the initial mesh. Upwind numerical flux is used at the elements' faces. $\mu = 0.001$ is used in the augmented Lagrangian for backtracking, the weight on mesh distortion is set to $\alpha = 5\times10^{-6}$ and the initial regularization in Hessian is set to $\varepsilon_{s0} = 1.0$ and $\varepsilon_{u0} = 0$ in algorithm \ref{algo_eq}. The control variables, $\vect{s}$, are highlighted in Fig. \ref{control_nodes_cshock}. Boundary nodes are constrained to slide along the boundary. The node at $(0.375, 1)$ is held fixed to ensure that the boundary condition on $u_0$ is evaluated accurately.

\begin{figure}[H]
    \centering
    \includegraphics[scale=0.4]{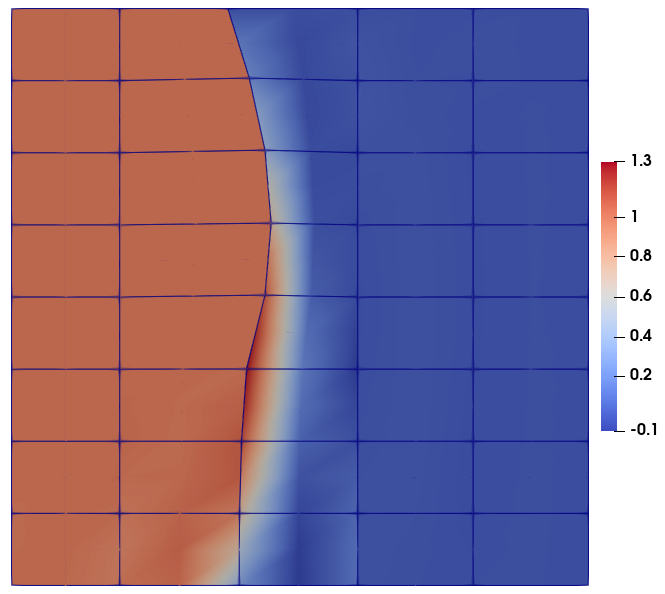}
    \caption{Optimal mesh and solution $u_0$ obtained using adjoint-based full space shock tracking on a grid of order 1.}
    \label{optimal_mesh_q1_cshock_goaloriented}
\end{figure}

\begin{figure}[H]
    \centering
    \begin{subfigure}{0.55\textwidth}
    \includegraphics[width=1\linewidth]{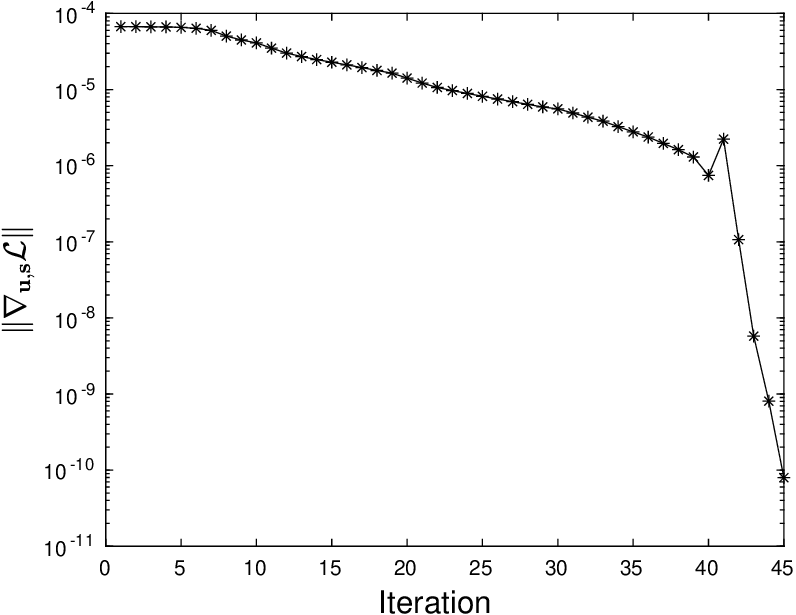}
    \caption{}
    \end{subfigure}
    
    \begin{subfigure}{0.55\textwidth}
    \includegraphics[width=1\linewidth]{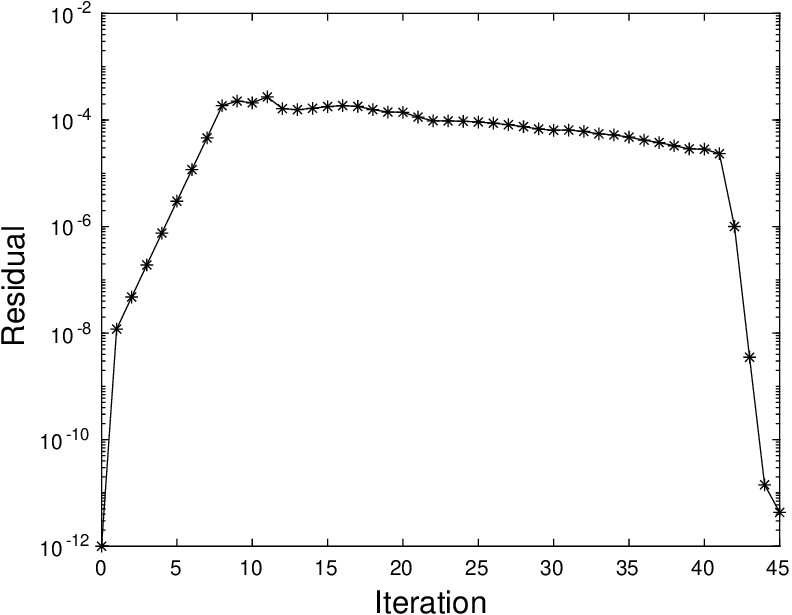}
    \caption{}
    \end{subfigure}
    \caption{Convergence of Lagrangian gradient, $\vect{\nabla \mathcal{L}}$, for goal-oriented shock tracking with a grid of order 1. (a) convergence of $\|\vect{\nabla_{u,s} \mathcal{L}} \|$. (b) convergence of $\|\vect{\nabla_{\lambda}\mathcal{L}}\| =\|\vect{r(u,x)}\|$.}
    \label{gradient_q1_goaloriented_cshock}
\end{figure}

\begin{figure}[H]
    \centering
    \includegraphics[scale=0.4]{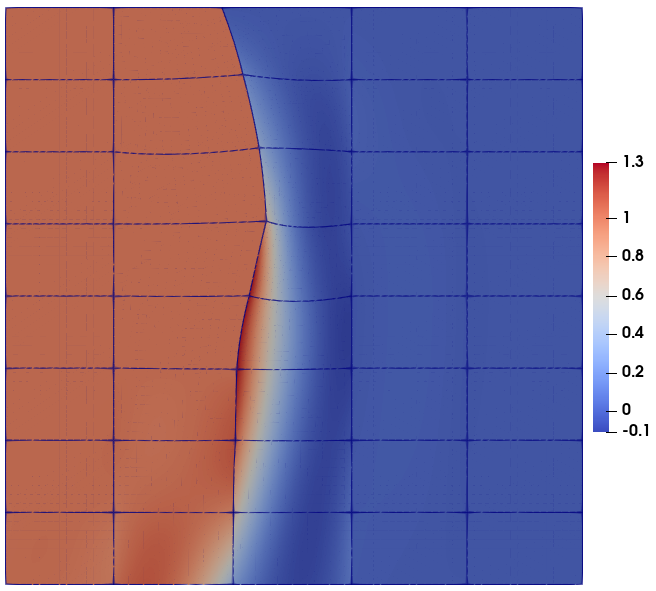}
    \caption{Optimal mesh and solution $u_0$ obtained using adjoint-based full space shock tracking on a grid of order 2.}
    \label{optimal_mesh_q2_cshock_goaloriented}
\end{figure}

\begin{figure}[H]
    \centering
    \begin{subfigure}{0.55\textwidth}
    \includegraphics[width=1\linewidth]{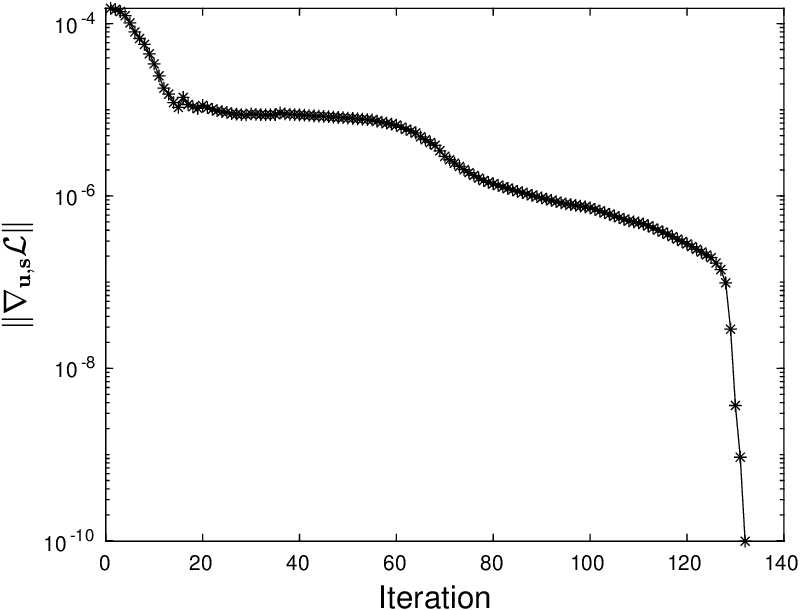}
    \caption{}
    \end{subfigure}
    
    \begin{subfigure}{0.55\textwidth}
    \includegraphics[width=1\linewidth]{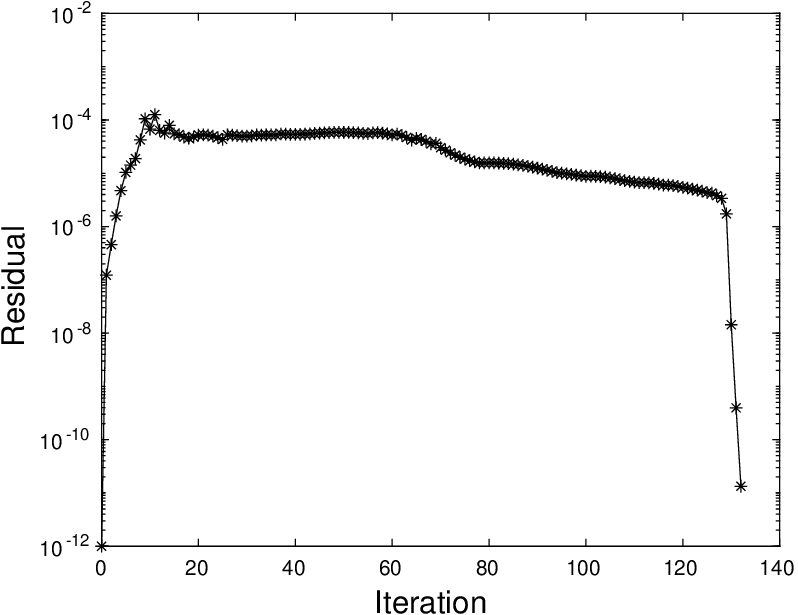}
    \caption{}
    \end{subfigure}
    \caption{Convergence of Lagrangian gradient, $\vect{\nabla \mathcal{L}}$, for goal-oriented shock tracking with a grid of order 2. (a) convergence of $\|\vect{\nabla_{u,s} \mathcal{L}} \|$. (b) convergence of $\|\vect{\nabla_{\lambda}\mathcal{L}}\| =\|\vect{r(u,x)}\|$.}
    \label{gradient_q2_goaloriented_cshock}
\end{figure}

Figs. \ref{optimal_mesh_q1_cshock_goaloriented} and \ref{gradient_q1_goaloriented_cshock} show respectively the optimal mesh and the convergence of the Lagrangian gradient for the mesh of grid degree 1, while Figs. \ref{optimal_mesh_q2_cshock_goaloriented} and \ref{gradient_q2_goaloriented_cshock} show the same for the mesh of grid degree 2. As seen from Figs. \ref{optimal_mesh_q1_cshock_goaloriented} and \ref{optimal_mesh_q2_cshock_goaloriented}, the optimizer aligns elements' faces with the shock in the top-half of the domain.

%%%%%%%%%%%%
%%%%%%%%%%%%

%%%%%%%%%%%%
\begin{figure}[H]
    \centering \includegraphics[scale=0.6]{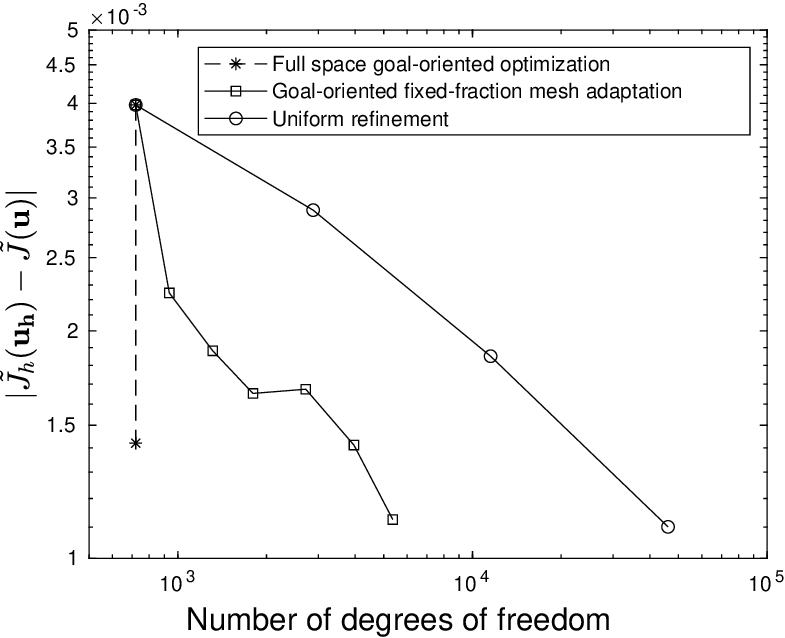}
\caption{Convergence of the functional error vs degrees of freedom for various mesh adaptation strategies.}
\label{functional_error_cshock_mesh_adaptations}
\end{figure}

\begin{figure}[H]
    \centering
    \includegraphics[scale=0.4]{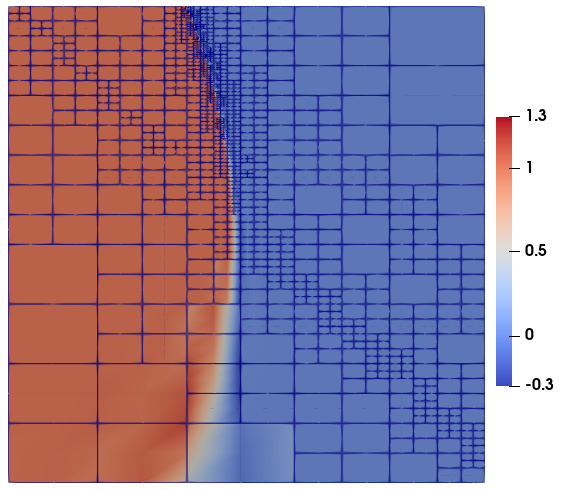}
    \caption{Mesh obtained from fixed-fraction goal-oriented $h$-adaptation.}
    
    \label{fixedfraction_cshock}
\end{figure}

%%%%%%%%%%%%
In order to compute the error in the functional, note that the exact solution $u_1$ can be obtained analytically and is:
\begin{equation}
\label{u1_exact}
    u_1(x,y) = \begin{cases}
        1.3 - y, & y>1-x, \; x \leq -0.4y^2 + 0.4y + 0.375. \\
        0.3 + x, & y \leq 1-x, \;x \leq -0.4y^2 + 0.4y + 0.375. \\
        0.3 + \tilde{x}, & y\leq 1-x, \; x > -0.4y^2 + 0.4y + 0.375. \\ 
        1.3 - \tilde{y}, & x > -0.4y^2 + 0.4y + 0.375, \; 1-x<y< 1.375-x. \\
        0.3, & y \geq 1.375-x.
    \end{cases}
\end{equation}
where
\begin{equation}
    \tilde{y} = \frac{1.4 - \sqrt{1.4^2 + 1.6(0.375-x-y)}}{0.8}, \quad \tilde{x} = -0.4\tilde{y}^2 + 0.4\tilde{y} + 0.375.
\end{equation}
It can be verified that the exact $u_1$ in Eq. \ref{u1_exact} satisfies Eqs. \ref{u1_eq_boundary} and \ref{u1_eq_volume}. Since the exact $u_1$ is irrational and the functional cannot be integrated analytically, the true value of the functional, $\tilde{J}(\vect{u})$, is approximated using Gauss-Legendre quadrature rule in Eq. \ref{modified_functional_testcase2} with 500 quadrature points along each element's face on $\Gamma_{\text{right}}$. 
Fig. \ref{functional_error_cshock_mesh_adaptations} shows the reduction in the functional error vs the number of degrees of freedom for goal-oriented optimization-based $r$-adaptation with a grid of degree 2, goal-oriented fixed-fraction $h$-adaptation with $10\%$ refinement fraction and uniform $h$-refinement. As seen from Fig. \ref{functional_error_cshock_mesh_adaptations}, fixed-fraction and uniform refinement methods require around 7 and 64 times more degrees of freedom respectively to achieve functional errors similar to the optimization-based approach. Fig. \ref{fixedfraction_cshock} shows the mesh obtained from fixed-fraction adaptation. 

% Include optimiation vs fixed-fraction vs uniform refinement plot.

\begin{figure}[H]
    \centering
    \includegraphics[scale=0.4]{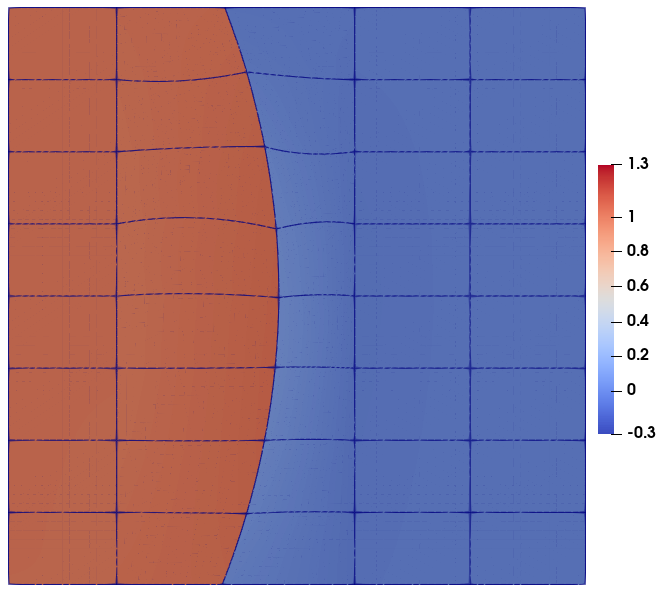}
    \caption{Optimal mesh and solution $u_0$ obtained using residual-based full space shock tracking on a grid of order 2.}
    \label{optimal_mesh_q2_cshock_residualbased}
\end{figure}

\begin{figure}[H]
    \centering
    \begin{subfigure}{0.55\textwidth}
    \includegraphics[width=1\linewidth]{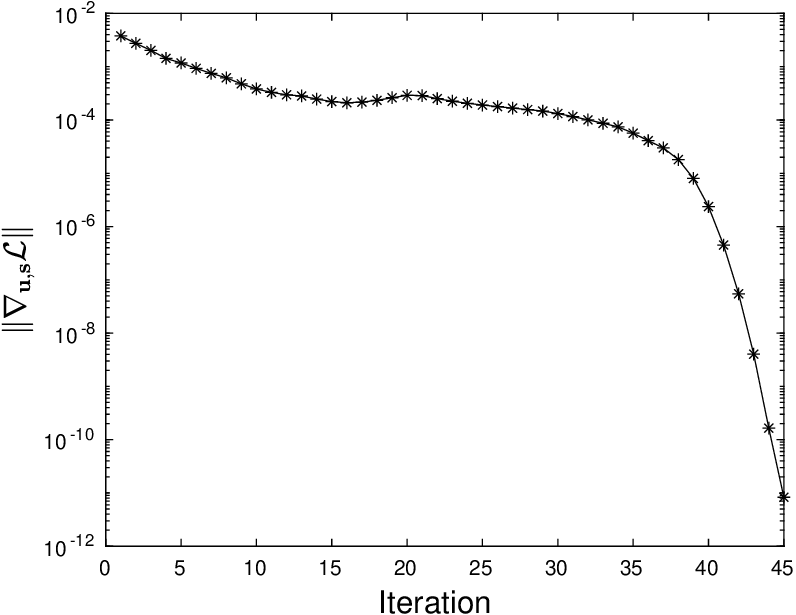}
    \caption{}
    \end{subfigure}
    
    \begin{subfigure}{0.55\textwidth}
    \includegraphics[width=1\linewidth]{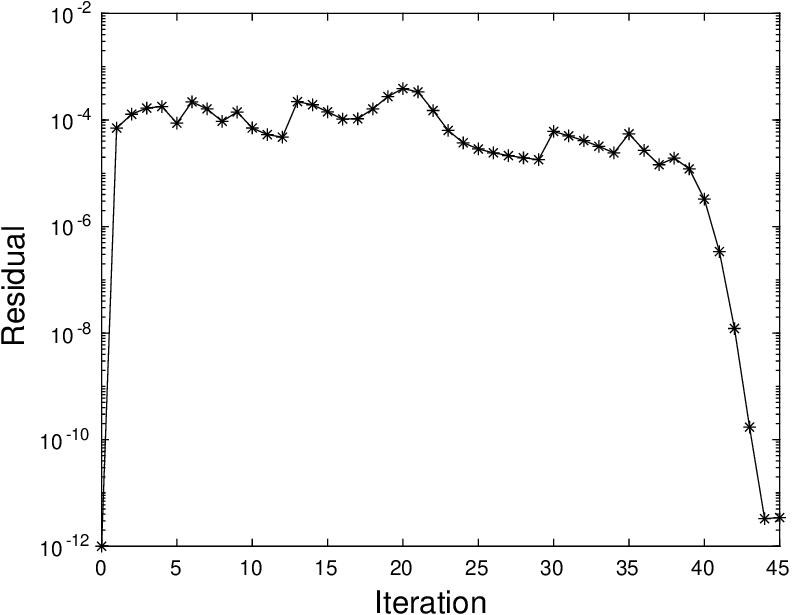}
    \caption{}
    \end{subfigure}
    \caption{Convergence of Lagrangian gradient, $\vect{\nabla \mathcal{L}}$, for residual-based shock tracking with a grid of order 2. (a) convergence of $\|\vect{\nabla_{u,s} \mathcal{L}} \|$. (b) convergence of $\|\vect{\nabla_{\lambda}\mathcal{L}}\| =\|\vect{r(u,x)}\|$.}
    \label{gradient_q2_residualbased_cshock}
\end{figure}

We also run the residual-based shock tracking of \cite{zahr2} for a grid of order 2, with smooth approximation of the upwind flux as proposed in \cite{zahr2} and for $p=1$. Fig. \ref{optimal_mesh_q2_cshock_residualbased} shows the optimal mesh and Fig. \ref{gradient_q2_residualbased_cshock} shows the convergence of the Lagrangian gradient. As seen from Figs. \ref{optimal_mesh_q2_cshock_goaloriented} and \ref{optimal_mesh_q2_cshock_residualbased}, the residual-based optimization tracks the entire shock while the goal-oriented optimization can track the shock in the region important for evaluating the functional.

%%%%%%%%%%%%

\subsection{Inviscid supersonic flow over a cylinder}
Finally, full space goal-oriented shock tracking is used to track the bow shock forward of a cylinder placed in supersonic inviscid flow. Tracking the entire shock with the optimizer results in the convergence of the DG residual without any addition of artificial dissipation. Since the governing Euler PDE is solved numerically without the use of any artificial means to stabilize the shock, the method yields optimal convergence orders of $\mathcal{O}(h^{p+1})$.

The computational domain is parameterized similar to the test case in Shu \cite{Shu1998}:
\begin{equation}
\begin{split}
    &  x = -[3 - 2\tilde{\xi}]\cos(\theta(2\tilde{\eta}-1)) \\
    &  y = [6 - 5\tilde{\xi}]\sin(\theta(2\tilde{\eta}-1)) 
\end{split}
\end{equation}
with $\theta = 5\pi/12$ and $(\tilde{\xi},\tilde{\eta}) \in [0\;\;1]^2$. It consists of a cylinder of radius 1 (at $\tilde{\xi}=1$) with slip-wall boundary condition, supersonic inflow at the left boundary ($\tilde{\xi}=0$) and supersonic outflow at the top and bottom boundaries ($\tilde{\eta}=0$ and $\tilde{\eta}=1$). The flow is governed by the Euler equations:
\begin{equation}
\label{eq_euler}
    \vect{\nabla} \cdot \vect{F(u)} = 0 \quad \text{on} \quad \Omega.
\end{equation}
with the conservative solution $\vect{u}$ and the convective flux $\vect{F} = [\vect{F_0} \quad \vect{F_1}]$:
\begin{equation}
    \vect{u} = \begin{bmatrix}
        \rho \\
        \rho u \\
        \rho v \\
        \rho E
    \end{bmatrix},
    \quad
    \vect{F_0} = \begin{bmatrix}
        \rho u \\
        \rho u^2 + p \\
        \rho uv \\
        \rho H u
    \end{bmatrix},
    \quad
    \vect{F_1} = \begin{bmatrix}
        \rho v \\
        \rho uv \\
        \rho v^2 + p \\   
        \rho H v
    \end{bmatrix}.
\end{equation}
$\rho$ is the density, $\vect{v} = [u\quad v]^T$ is the velocity vector, $E$ is the specific total energy, $p = (\gamma-1)\left(\rho E -\frac{1}{2}\rho \|\vect{v}\|^2 \right)$ is the pressure, $\gamma = 1.4$ is the specific heat ratio of air and $H = (\rho E + p)/\rho$ is the enthalpy. Freestream conditions similar to Corrigan et al. \cite{mdg1} are used, where the flow is moving towards the cylinder at an angle of attack of $0^0$ with $M_{\infty}=3$, $p_{\infty} = 1$ and $\rho_{\infty} = 1.4$. Since the flow is inviscid, enthalpy is conserved and is expected to be $H_{\text{exact}} = H_{\infty} = 7$ throughout the domain. For this test case, we choose the functional as the mass flow rate out of the domain:
\begin{equation}
    J(u) = \int_{\Gamma_{outflow}} \rho \vect{v}\cdot \vect{n} d\Gamma.
\end{equation}

DG discretization described in section \ref{DG} is used to discretize Eq. \ref{eq_euler}. We use the adjoint consistent DG for Euler equations \cite{hartmann_adjoint_consistent}, as it has been noted in \cite{hartmann_bali, hartmann_adjoint_consistent} that an adjoint inconsistent discretization results in irregularities in the adjoint solution. To that end, at the surface of the cylinder, solution using the slip wall boundary condition, $\vect{u}_{\Gamma}(\vect{u})$, is obtained from $\vect{u}$ by replacing the velocity in the interior, $\vect{v^+}$, with $\vect{v}_{\Gamma} = \vect{v^+} - (\vect{v^+}\cdot \vect{n})\vect{n}$. Moreover, the functional is evaluated using $\vect{u}_{\Gamma}$. Roe flux with entropy fix of Hartmann and Leicht \cite{hartmann_leicht} is used at the faces, with the Roe-averaged eigenvalues ($\lambda_r^i,\;i=1,2,3,4$) modified as:
\begin{equation}
    \tilde{\lambda_r^i} = \begin{cases}
        & \frac{\left(\lambda_r^i\right)^2 + \delta^2}{2\delta}, \quad |\lambda_r^i| < \delta \\
        & |\lambda_r^i|, \quad \text{otherwise} 
    \end{cases}
\end{equation}
with $\delta = 0.1 \max\limits_{i} \{|\lambda_r^i|\}$. Since this modification does not satisfy the Rankine-Hugoniot jump condition at the shock if $|\lambda_r^i| < \delta$, we use pure upwinding at the face on the shock as described later. Note that an entropy fix is necessary for this test case since the standard Roe flux is known to not be accurate near sonic rarefaction \cite{toro_riemannsolvers} and can cause carbuncle instability. Further, we observed that standard Roe and HLLC fluxes can fail to impose upwinding at the shock face when the mesh is aligned with the shock as they solve the Riemann problem using approximations of the wave speeds. The same was also observed in Huang et al.\cite{zahr7}. 

Careful treatment of the wall boundary is essential to obtain orders of convergence for this test case. It was found that using either the exact normals and surface Jacobians \cite{exact_normals} or superparametric grids \cite{PhilipZ} yields optimal convergence orders for cases involving curved boundaries. In the current work, we use the former and, similar to \cite{mdg1}, the physical quadrature on the cylinder obtained from the polynomial representation of the grid ($\vect{x(\xi)}$) is projected onto the cylinder:
\begin{equation}
    \vect{P(x)} = \frac{\vect{x(\xi)}}{\|\vect{x(\xi)}\|}.
\end{equation}
The Jacobian can be computed using chain rule, $\vect{J} = \frac{\partial \vect{P(x)}}{\partial \vect{x}} \frac{\partial \vect{x}}{\partial \vect{\xi}}$, and the cofactor of this Jacobian is used to evaluate the normal and the surface integral in Eq. \ref{weak_form_final}.

\begin{figure}[H]
    \centering
    \includegraphics[scale=0.3]{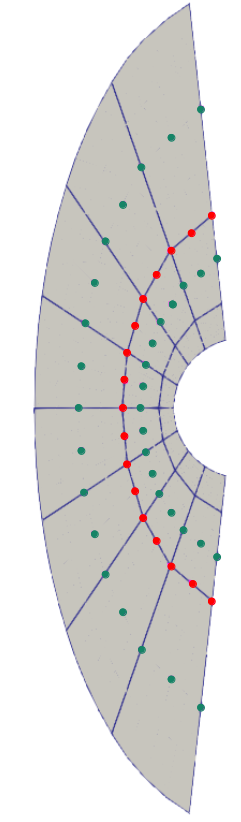}
    \caption{Control nodes used in optimization are highlighted in red. Nodes in green are displaced such that they are at the center of their respective elements/faces. }
    \label{cylinder_controlnodes}
\end{figure}

Figure \ref{cylinder_controlnodes} shows the initial mesh and the control nodes used for optimization. The control nodes are constrained to move radially i.e. along the line joining the node on the cylinder to the corresponding node on the inflow boundary so that the shock tracking problem has a unique solution.   

\begin{figure}[H]
    \centering
    \includegraphics[scale=0.4,trim={8cm 0 18cm 0},clip]{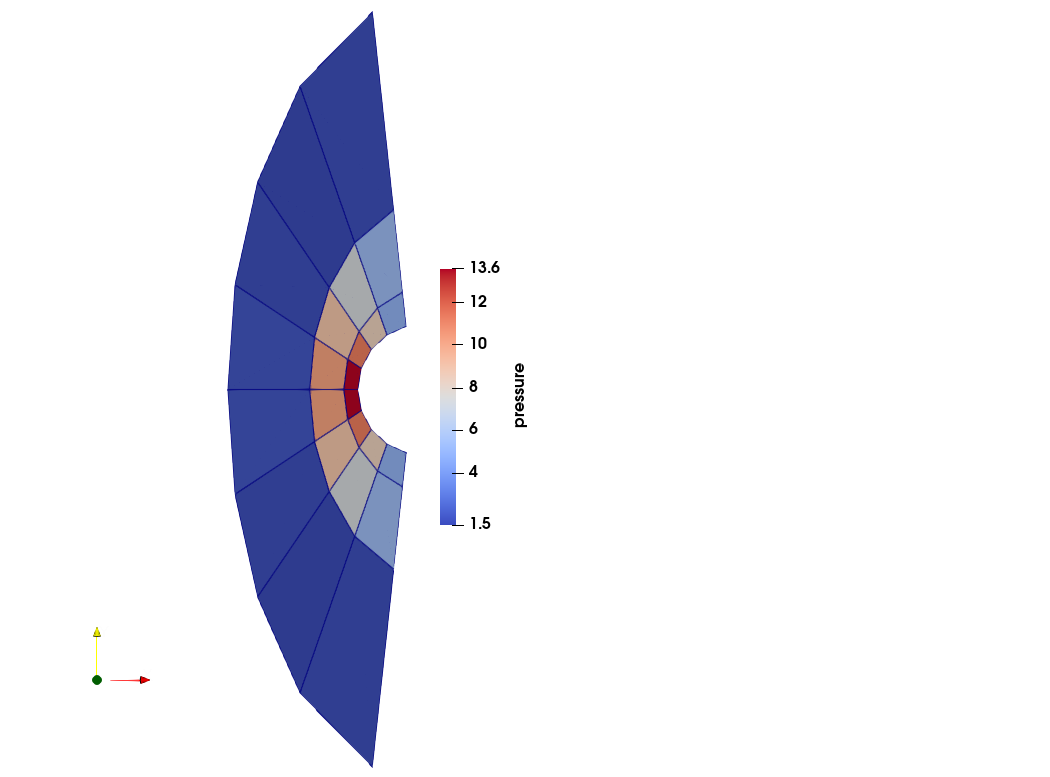}
    \caption{Initial mesh with the pressure field using $p=0$ solution on $q=1$ grid.}
    \label{cylinder_initial_mesh}
\end{figure}

In order to initialize the shock tracking problem, we used a $p=0$ solution converged on a coarse grid of order $q=1$, containing 24 cells as shown in Figure \ref{cylinder_initial_mesh}. Next, the solution was interpolated to $p=1$, keeping $q=1$. The full space optimizer was run for 4 iterations using Roe with entropy fix throughout the domain. Finally, the solution was interpolated to a grid of order $q=2$ while keeping $p=1$. At this stage, the DG residual used upwinding at the faces containing control nodes while using Roe with entropy fix elsewhere, and the full space optimizer was run until convergence to a tolerance of $10^{-10}$. 

A penalty of $\mu = 0.001$ was used in the augmented Lagrangian for backtracking. Initial regularization parameters of control ($\varepsilon_{s0}$) and simulation ($\varepsilon_{u0}$) were both set to 1. Regularization matrix $\vect{D}$ for control was chosen to be the Poisson equation's stiffness matrix for this test case, as discussed in section \ref{section_regularization}. Weight on mesh distortion ($\alpha$) was set to 0. Due to the size of the test case, in the Hessian on the left-hand side of the KKT system in Eq. \ref{regularized_system}, we use the preconditioner of $\vect{R_u}$ while evaluating vector products with the derivatives of the adjoint in Eq. \ref{adjoint_derivatives}. However, we still use the exact $\vect{R_u}$ in adjoint derivatives while forming the Lagrangian gradient on the right-hand side of Eq. \ref{regularized_system}.

For this test case, regularizing the solution was important to reduce oscillations arising from the Gibbs phenomenon as the optimizer converges. The basic idea is to have higher regularization on the flow solution as the mesh moves and update the regularization parameter as described in section \ref{section_regularization} such that the regularization on the flow solution approaches zero when the mesh has almost tracked the shock. After the full space optimization converged on the coarse grid for $p=1$ and $q=2$, the mesh was isotropically refined and the solution was interpolated to the fine mesh before running the optimizer again. For $p=1$, $q=2$, the optimization was run on four sets of grids to test the orders of convergence: grid 1 (24 cells), grid 2 (96 cells), grid 3 (384 cells), and grid 4 (1536 cells).

Next, using the converged optimal mesh and solution from various grids of $(p=1, q=2)$, we interpolate the solution to $p=2$ and re-run the full space optimizer to obtain the optimal mesh for $(p=2, q=2)$. For $(p=2, q=2)$, optimization is performed on grids 1 to 3. 

\begin{figure}[H]
    \centering
    \begin{subfigure}[b]{0.18\textwidth}
    \centering
    \includegraphics[width=\textwidth,trim={8cm 0 22cm 0},clip]{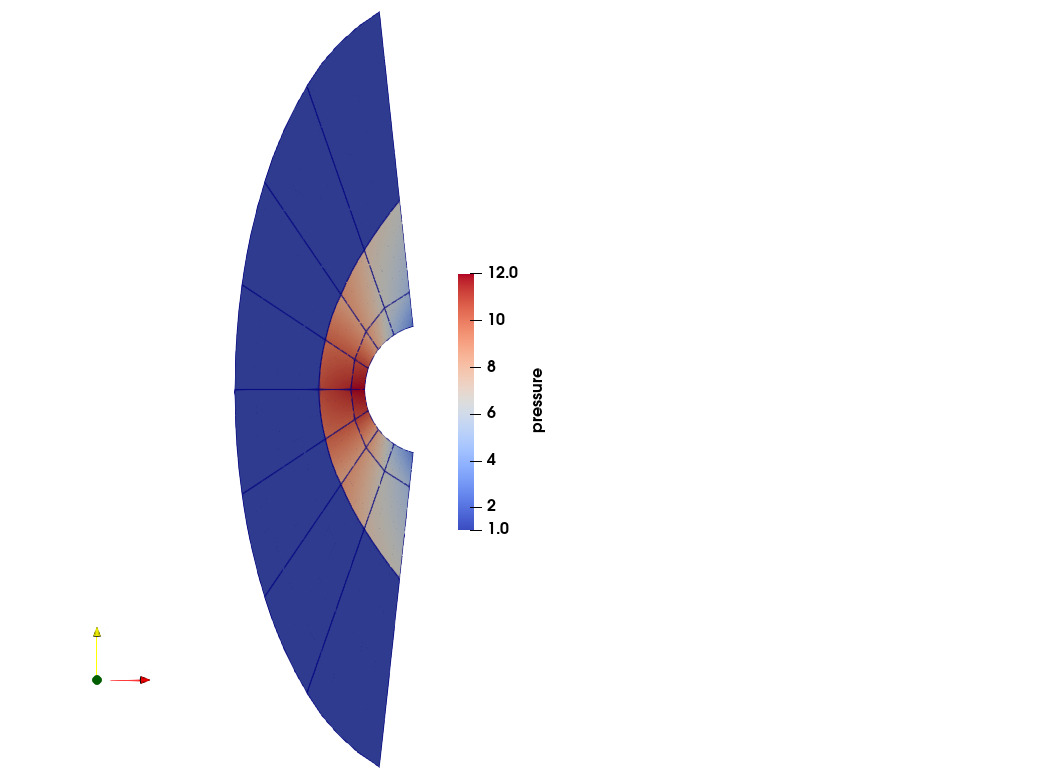}
    \caption{Grid 1}
    \label{p1_grid1}      
    \end{subfigure}    
    \hfill
    \begin{subfigure}[b]{0.18\textwidth}
    \centering
    \includegraphics[width=\textwidth,trim={8cm 0 22cm 0},clip]{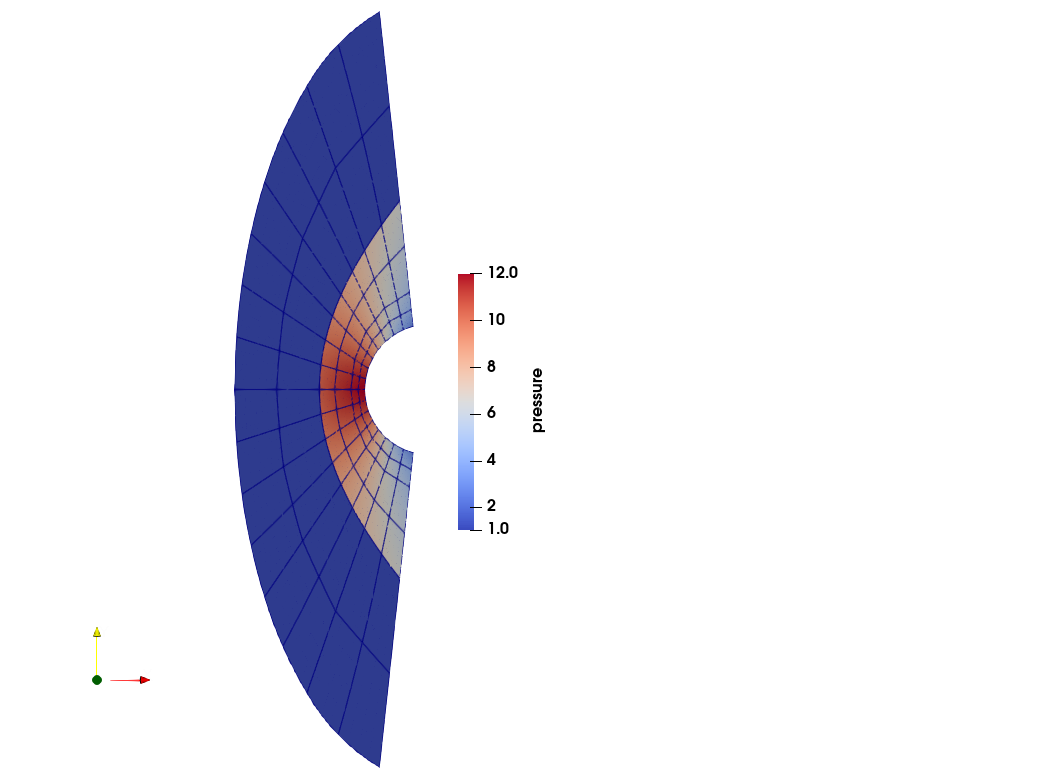}
    \caption{Grid 2}
    \label{p1_grid2}      
    \end{subfigure}
        \hfill
    \begin{subfigure}[b]{0.18\textwidth}
    \centering
    \includegraphics[width=\textwidth,trim={8cm 0 22cm 0},clip]{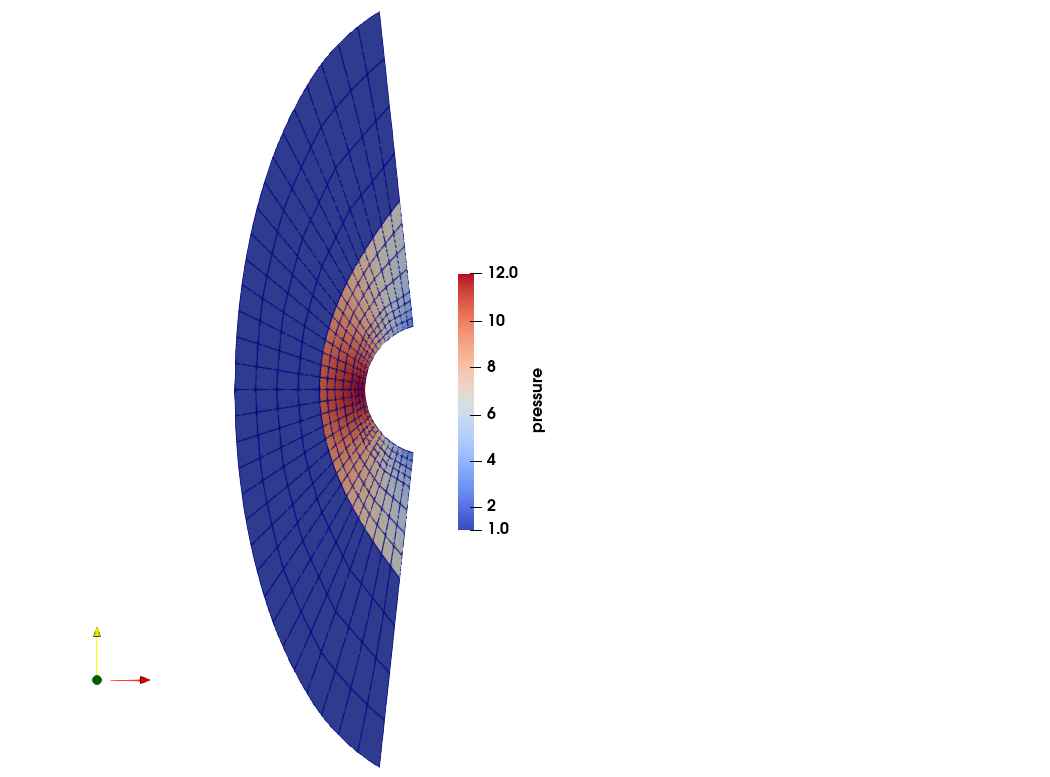}
    \caption{Grid 3}
    \label{p1_grid3}      
    \end{subfigure}
            \hfill
    \begin{subfigure}[b]{0.18\textwidth}
    \centering
    \includegraphics[width=\textwidth,trim={8cm 0 22cm 0},clip]{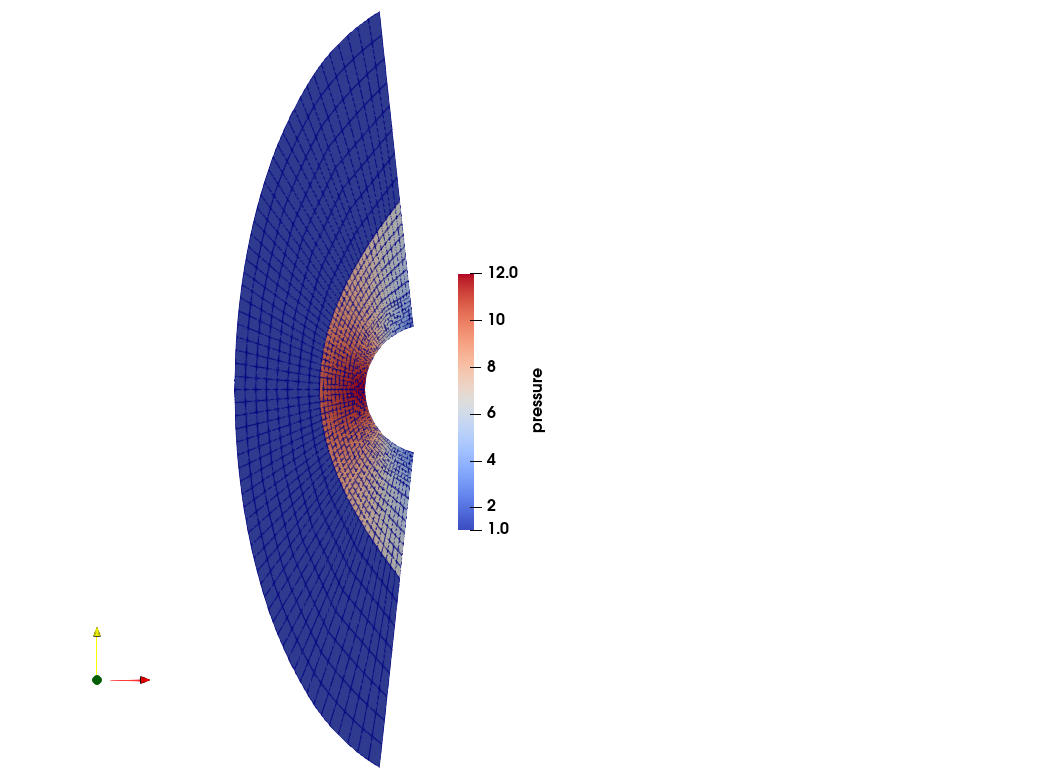}
    \caption{Grid 4}
    \label{p1_grid4}      
    \end{subfigure}
    \hfill
    \begin{subfigure}[b]{0.14\textwidth}
    \centering
    \includegraphics[width=\textwidth,trim={15cm 0 17cm 0},clip]{p1_mesh_4.png}      
    \end{subfigure}

    \caption{Optimal meshes with the pressure field on a sequence of grids for $(p=1, q=2)$.}
    \label{p1q2_cylinder_meshes}
\end{figure}

\begin{figure}[H]
    \centering
    \begin{subfigure}[b]{0.18\textwidth}
    \centering
    \includegraphics[width=\textwidth,trim={8cm 0 22cm 0},clip]{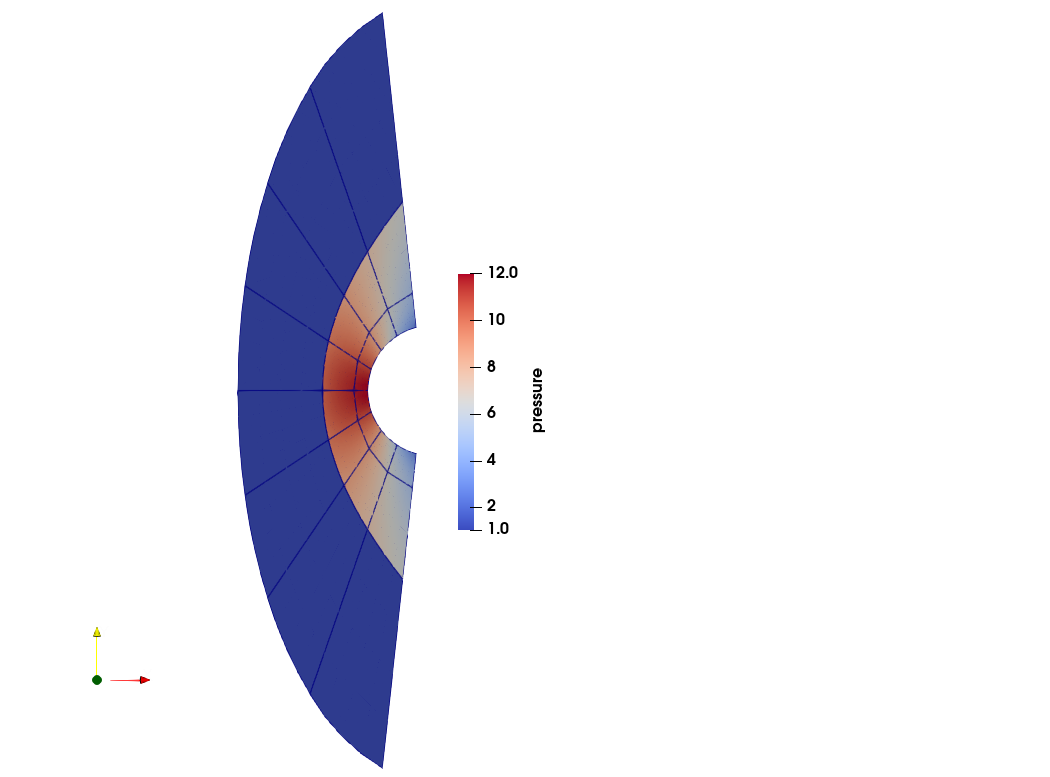}
    \caption{Grid 1}
    \label{p2_grid1}      
    \end{subfigure}    
    \hfill
    \begin{subfigure}[b]{0.18\textwidth}
    \centering
    \includegraphics[width=\textwidth,trim={8cm 0 22cm 0},clip]{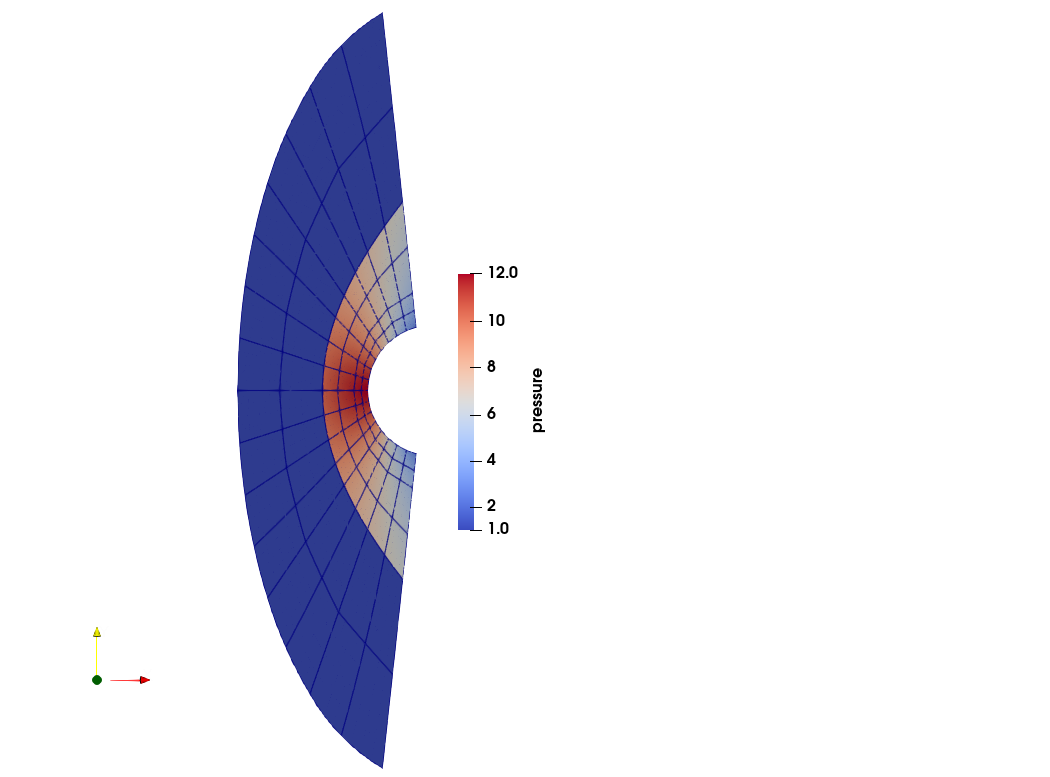}
    \caption{Grid 2}
    \label{p2_grid2}      
    \end{subfigure}
        \hfill
    \begin{subfigure}[b]{0.18\textwidth}
    \centering
    \includegraphics[width=\textwidth,trim={8cm 0 22cm 0},clip]{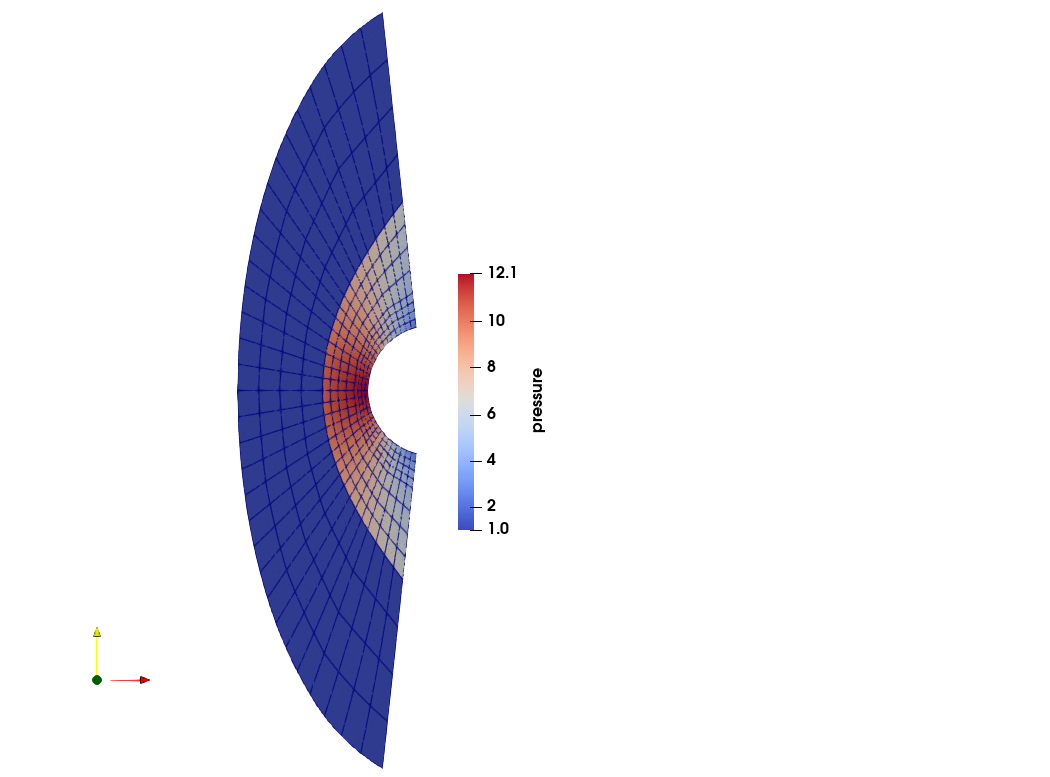}
    \caption{Grid 3}
    \label{p2_grid3}      
    \end{subfigure}
            \hfill
    \begin{subfigure}[b]{0.14\textwidth}
    \raggedleft
    \includegraphics[width=\textwidth,trim={15cm 0 17cm 0},clip]{p2_mesh_3.png}      
    \end{subfigure}

    \caption{Optimal meshes with the pressure field on a sequence of grids for $(p=2, q=2)$.}
    \label{p2q2_cylinder_meshes}
\end{figure}

\begin{figure}[H]
    \centering
    \begin{subfigure}[b]{0.32\textwidth}
    \centering
    \includegraphics[width=\textwidth,trim={8cm 0 17cm 0},clip]{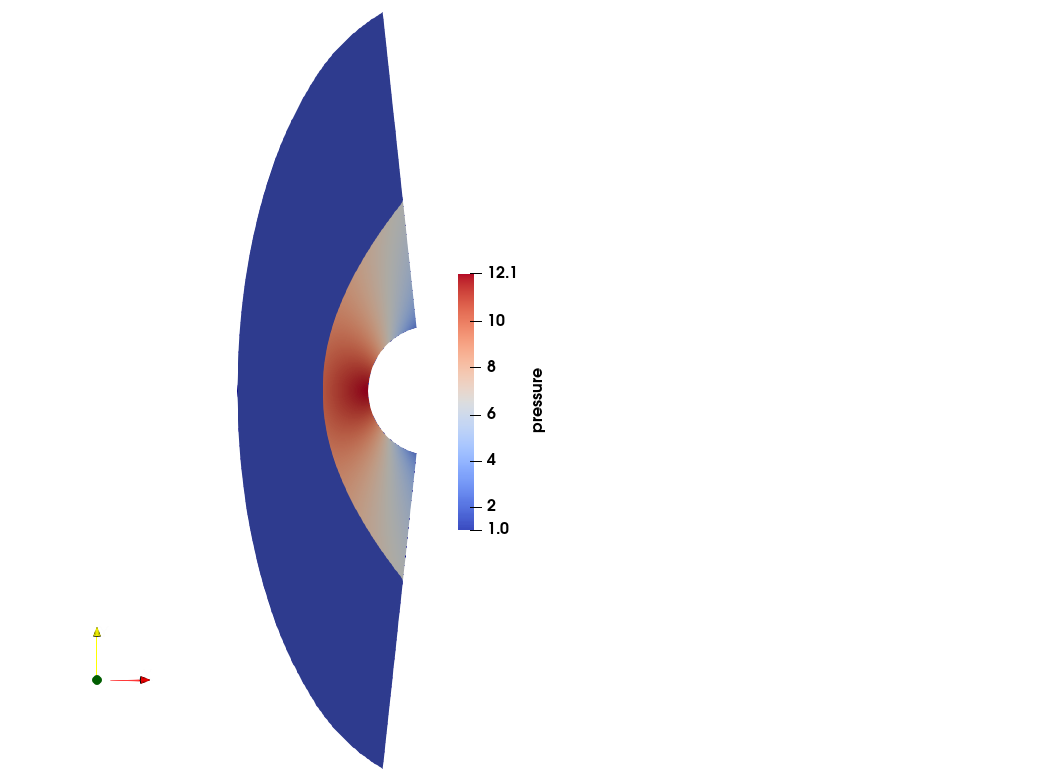}
    \caption{Pressure}
    \label{p2_pressure_field}      
    \end{subfigure}    
    \hfill
    \begin{subfigure}[b]{0.35\textwidth}
    \centering
    \includegraphics[width=\textwidth,trim={7cm 0 17cm 0},clip]{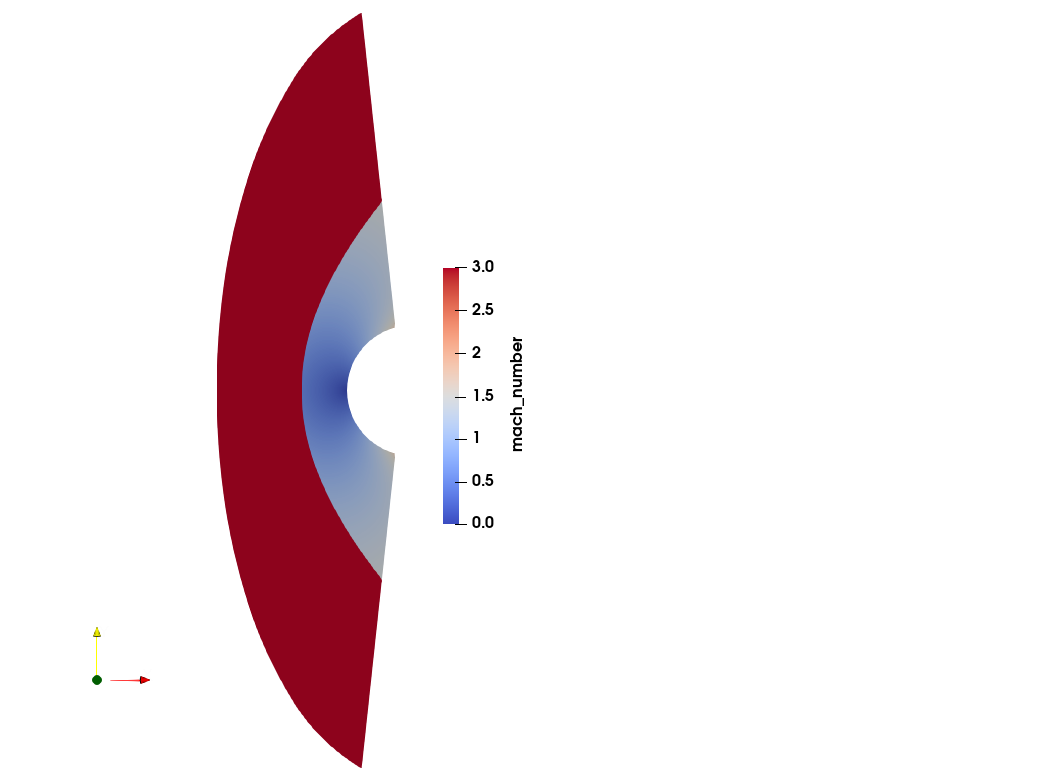}
    \caption{Mach number}
    \label{p2_mach_field}      
    \end{subfigure}

    \caption{Plot of pressure and Mach number on the optimal grid 3 of $(p=2, q=2)$.}
    \label{cylinder_pressure_mach_field}
\end{figure}

Figure \ref{p1q2_cylinder_meshes} shows the optimal mesh and pressure field on various grids for $(p=1, q=2)$ while Figure \ref{p2q2_cylinder_meshes} shows the same for $(p=2, q=2)$. Figure \ref{cylinder_pressure_mach_field} shows the plot of the pressure field and the Mach number without the mesh. The mesh tracks the shock in about 30 full space iterations for all cases. 

\begin{figure}[H]
    \centering
    \begin{subfigure}{0.55\textwidth}
    \includegraphics[width=1\linewidth]{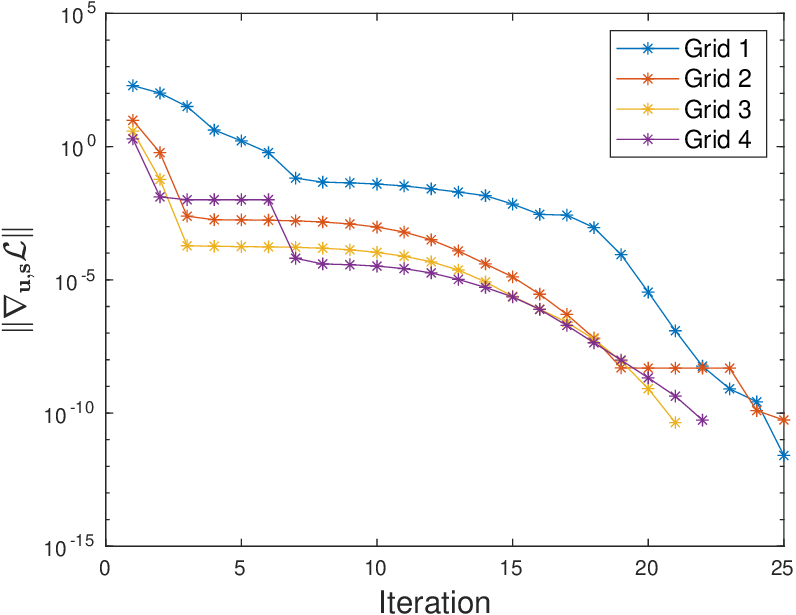}
    \caption{}
    \end{subfigure}
    
    \begin{subfigure}{0.55\textwidth}
    \includegraphics[width=1\linewidth]{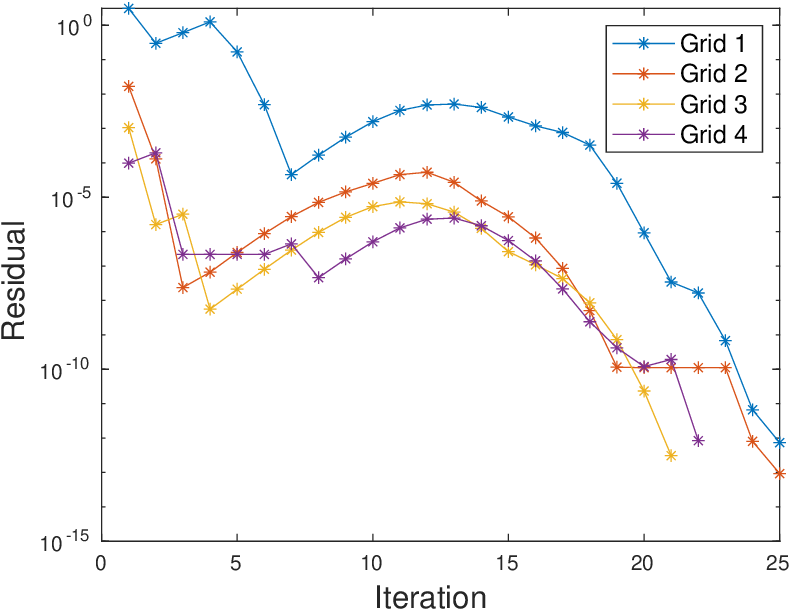}
    \caption{}
    \end{subfigure}
    \caption{Convergence of Lagrangian gradient, $\vect{\nabla \mathcal{L}}$, for goal-oriented shock tracking with $(p=1, q=2)$. (a) convergence of $\|\vect{\nabla_{u,s} \mathcal{L}} \|$. (b) convergence of $\|\vect{\nabla_{\lambda}\mathcal{L}}\| =\|\vect{r(u,x)}\|$.}
    \label{gradient_p1q2_cylinder}
\end{figure}

\begin{figure}[H]
    \centering
    \begin{subfigure}{0.55\textwidth}
    \includegraphics[width=1\linewidth]{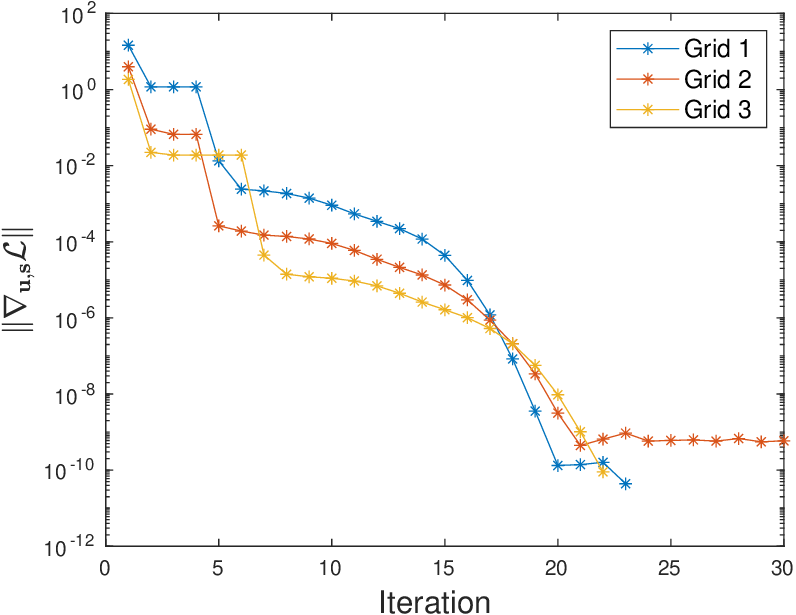}
    \caption{}
    \end{subfigure}
    
    \begin{subfigure}{0.55\textwidth}
    \includegraphics[width=1\linewidth]{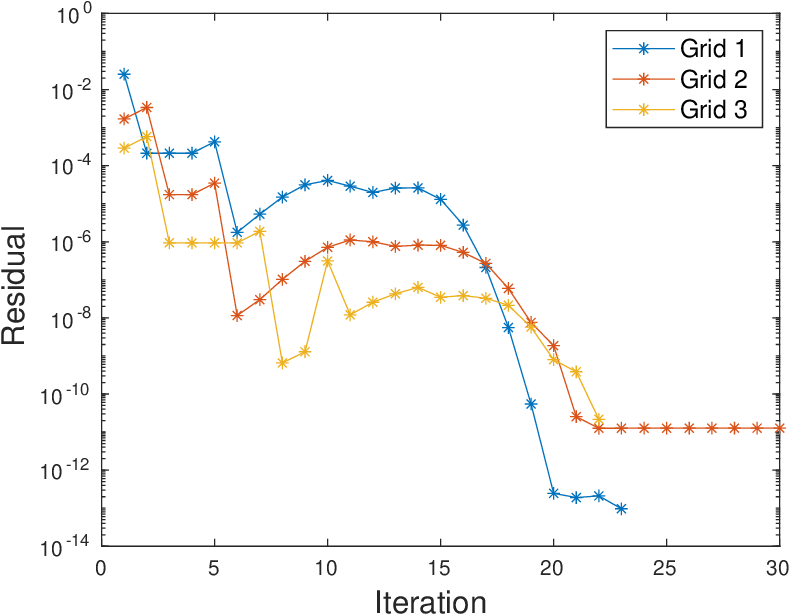}
    \caption{}
    \end{subfigure}
    \caption{Convergence of Lagrangian gradient, $\vect{\nabla \mathcal{L}}$, for goal-oriented shock tracking with $(p=2, q=2)$. (a) convergence of $\|\vect{\nabla_{u,s} \mathcal{L}} \|$. (b) convergence of $\|\vect{\nabla_{\lambda}\mathcal{L}}\| =\|\vect{r(u,x)}\|$.}
    \label{gradient_p2q2_cylinder}
\end{figure}

Figure \ref{gradient_p1q2_cylinder} shows the convergence of the gradient of the Lagrangian and the residual on the meshes used for $(p=1, q=2)$ case while Figure \ref{gradient_p2q2_cylinder} shows the same for $(p=2, q=2)$. DG residual converged below the tolerance of $10^{-10}$ for all cases while the Lagrangian gradient converged below $10^{-10}$ for all grids except grid 2 of $(p=2, q=2)$, where it stalled at $6\times10^{-10}$. We note that in contrast to the test case in \cite{mdg1} where the overdetermined residual cannot be solved exactly but can be minimized in a least squares sense, the current approach can achieve converged residual within the framework of full space optimization. We also experimented with the Gauss-Newton approximation of the Hessian and observed that it resulted in the Lagrangian gradient stalling at $\sim 10^{-8}$ or requiring a few more nonlinear iterations, as expected, for convergence compared to using second-order derivatives in the Hessian.

%%%%%%%%%%%%
%%%%%%%%%%%%%%%%%

%%%%%%%%%%%%
%%%%%%%%%%%%%%%%%

\begin{figure}[H]
    \centering
    \includegraphics[scale=0.65]{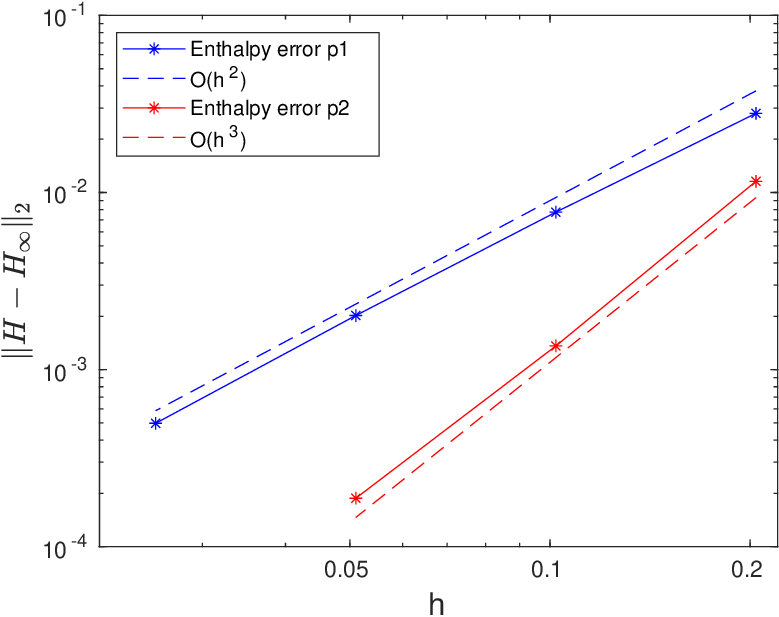}
    \caption{Convergence of the error in enthalpy using adjoint-based full-space shock tracking with respect to $h = 1/\sqrt{N_{\text{cells}}}$.}
    \label{cylinder_enthalpy}
\end{figure}

\begin{figure}[H]
    \centering
    \includegraphics[scale=0.65]{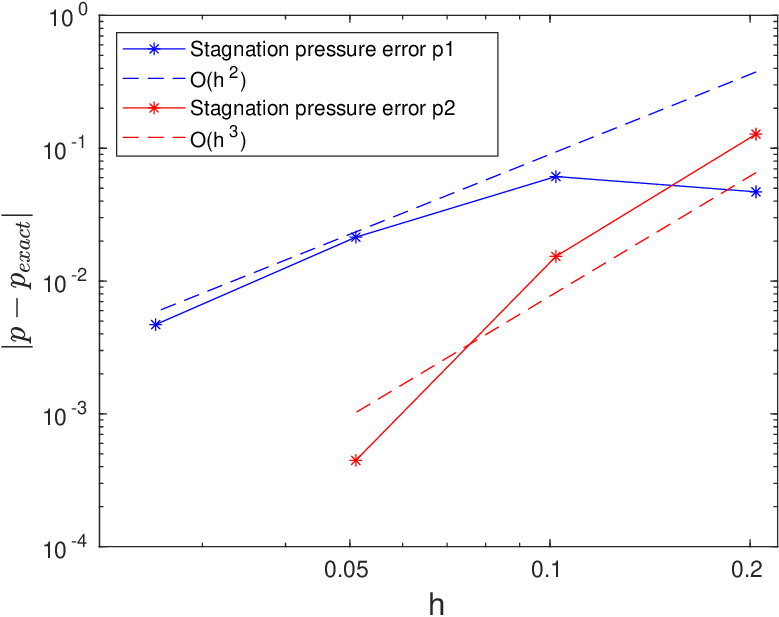}
    \caption{Convergence of the error in stagnation pressure using adjoint-based full space shock tracking with respect to $h = 1/\sqrt{N_{\text{cells}}}$.}
    \label{cylinder_pressure}
\end{figure}

Figure \ref{cylinder_enthalpy} shows the convergence of the $L_2$ norm of the error in enthalpy:
\begin{equation*}
    \|H(\vect{u_h}) - H_{\infty}\|_2 = \sqrt{\int_{\Omega} \left(H(\vect{u_h}) - H_{\infty}\right)^2 d\Omega},
\end{equation*}
while Figure \ref{cylinder_pressure} shows the convergence in the stagnation pressure. The exact stagnation pressure can be computed for this test case \cite{masatsuka} and is around 12.06096470126662. Expected orders of $\mathcal{O}(h^2)$ and $\mathcal{O}(h^3)$ are observed for both the error in enthalpy and stagnation pressure, with the convergence in enthalpy error being much smoother. This is in contrast to the commonly used artificial dissipation-based shock capturing approaches, which result in sub-optimal rates of convergence of errors.

%% file: 8Conclusion.tex
\section{Conclusion}
A framework for goal-oriented adjoint-based shock tracking was developed, using an objective function based on the adjoint-weighted residual. First and second-order derivatives of the objective function were derived and used in LNKS full space optimization to minimize functional error while converging the residual of the governing PDE simultaneously. Second-order derivatives were used to enhance the convergence of the optimizer. We observed that using second-order derivatives led to the optimizer converging at a lower tolerance with fewer iterations than the Gauss-Newton approximation, albeit at a higher cost. However, a formal analysis of cost needs to be performed to justify the use of second-order derivatives. The optimizer tracked shocks in regions with high adjoint-weighted residuals. The method was applied to several test cases involving convection and Euler equations to demonstrate functional-dependent shock tracking. Comparisons with goal-oriented fixed-fraction $h$-adaptation and uniform refinement revealed the potential of goal-oriented shock tracking to achieve lower functional errors on coarser meshes. We also noted the importance of using smooth adjoint solutions for the convergence of the optimizer. Tracking the shock without any added artificial dissipation recovered the optimal $\mathcal{O}(h^{p+1})$ rate of convergence of error in the presence of shocks. In contrast to the existing state-of-the-art shock tracking approaches \cite{mdg1,zahr2}, the current approach uses a functional-dependent objective function, second-order derivatives in the Hessian in place of the Gauss-Newton approximation, and full space preconditioners to improve scalability.  

While the proposed approach of goal-oriented shock tracking was able to reduce error in the functional on coarse meshes, we note that tracking the entire shock instead of a part of the shock would lead to similar results. However, the current work of using a goal-oriented error indicator as the objective function allows for a future extension of the approach beyond shock tracking, specifically to better resolve multiple flow features such as boundary layers that would be important for evaluating functionals. Several important issues need to be addressed in the future to make the proposed method a viable approach. Vector products with the derivatives of the adjoint, Eq. \ref{adjoint_derivatives}, require a linear solve with $\vect{R_u}$. While using a preconditioner of $\vect{R_u}$ for this purpose in the Hessian did not severely affect the convergence of the optimizer, a formal investigation into the approximation of $\vect{R_u}$ might be necessary. Using quadrilateral and hexahedral meshes improves the efficiency of computing DG residual \cite{tet_vs_hex} as the basis functions in higher-dimensions can be expressed as tensor-products of one-dimensional basis functions. However, mesh operations such as edge collapsing, edge splitting, and edge swapping can be non-trivial for such meshes. Including robust mesh operations in the optimizer will significantly enhance its capabilities to obtain optimal mesh, starting from an arbitrary initial mesh. While the test cases we presented in this work were limited to two-dimensions, derivatives of adjoint-weighted residual derived in this work and the framework of full space optimization are also applicable to three-dimensions. Convergence of the optimizer for three-dimensional test cases needs further investigation.  

%% file: 9Acknowledgements.tex
\section{Acknowledgements}
The first author would like to thank Alex Cicchino and Doug Shi-Dong for their helpful discussions. The authors would like to gratefully acknowledge the financial support from the Department of Mechanical Engineering, McGill University, and the Natural Sciences and Engineering Research Council of Canada (NSERC) Discovery Grant RGPIN-2019-04791. Computational resources of the Digital Research Alliance of Canada were used in this work.

\begin{comment}
\section{Declaration of competing interest}
The authors declare that they have no known competing ﬁnancial interests or personal relationships that could have
appeared to inﬂuence the work reported in this paper.
\end{comment}

%% file: 10Appendix.tex
\section{Continuous adjoint equation for the problem in section \ref{test_case_2}}
\label{continuous_adjoint_appendix}
This section contains the derivation of the adjoint equations, Eqs. \ref{eq_psi0_vol} - \ref{eq_psi1_boundary}, used in section \ref{test_case_2}. For the ease of notation, we denote $\vect{u} = [u_0 \quad u_1]^T$ and the primal PDE, Eqs. \ref{u0_eq_volume} and \ref{u1_eq_volume}, as
\begin{equation}
\label{appendix_pde}
    -\vect{\nabla}\cdot \vect{F} + \vect{s(u)} = \vect{0},
\end{equation}
with
\begin{equation}
\label{FS_substituted}
    \vect{F} = [\vect{F_0} \quad \vect{F_1}], \quad \vect{F_0} = \vect{\beta}u_0, \quad \vect{F_1} = \vect{\mathcal{V}}u_1, \quad \vect{s(u)} = \begin{bmatrix}
        0 \\
        u_0^2
    \end{bmatrix}.
\end{equation}
The functional in Eq. \ref{eq_functional_case2} can be expressed as
\begin{equation}
\label{appendix_functional}
    J(\vect{u}) = \int_{\Gamma} j(\vect{u}) d\Gamma, \quad j(\vect{u}) = \begin{cases}
        u_1^2, & \text{on}\; \Gamma_{\text{right}} \\
        0, & \text{on}\; \Gamma \setminus \Gamma_{\text{right}}
    \end{cases}.
\end{equation}
Using the adjoint solution $\vect{\psi} = [\psi_0 \quad \psi_1]^T$, the Lagrangian $\mathcal{L}$ is formed with Eqs. \ref{appendix_pde} and \ref{appendix_functional}:
\begin{equation}
    \mathcal{L} = \int_{\Gamma} j(\vect{u}) d\Gamma + \int_{\Omega} \vect{\psi^T}\left(-\vect{\nabla}\cdot \vect{F} + \vect{s(u)} \right)d\Omega.
\end{equation}
Integrating the term $\int_{\Omega}\vect{\psi^T}\left(\vect{\nabla}\cdot \vect{F}\right)d\Omega$ by-parts yields
\begin{equation}
    \mathcal{L} = \int_{\Gamma} \left(j(\vect{u}) -\vect{\psi^T} \left(\vect{F}\cdot \vect{n}\right) \right) d\Gamma + \int_{\Omega} \left(\vect{\nabla \psi} : \vect{F} + \vect{\psi^T}\vect{s(u)} \right) d\Omega.
\end{equation}

An arbitrary perturbation in the mesh causes perturbation $\vect{\delta u}$ in the solution and results in a perturbation in $\mathcal{L}$. We choose the adjoint $\vect{\psi}$ such that the variation is $\mathcal{L}$ due to $\vect{\delta u}$ is zero, i.e. $\vect{\psi}$ is such that 
\begin{equation}
\label{appendix_perturbation_eq}
    \int_{\Gamma} \left(\vect{j_u} - \vect{\psi^T} \left(\vect{F_u}\cdot \vect{n}\right) \right)\vect{\delta u} d\Gamma + \int_{\Omega} \left(\vect{\nabla \psi} : \vect{F_u} + \vect{\psi^T}\vect{s_u} \right)\vect{\delta u} d\Omega = 0.
\end{equation}
In the volume term of Eq.~\ref{appendix_perturbation_eq}, $\vect{\delta u}$ is arbitrary. Thus, we get
\begin{equation}
    \vect{\nabla \psi} : \vect{F_u} + \vect{\psi^T}\vect{s_u} = \vect{0}\quad \text{in} \quad \Omega,
\end{equation}
yielding Eqs. \ref{eq_psi0_vol} and \ref{eq_psi1_vol}:
\begin{equation}
    \begin{split}
    &\vect{\nabla \psi} : \vect{F}_{u_0} + \vect{\psi^T}\vect{s}_{u_0} = \vect{\nabla}\psi_0 \cdot \vect{\beta} + 2\psi_1u_0= 0, \\
    &\vect{\nabla \psi} : \vect{F}_{u_1} + \vect{\psi^T}\vect{s}_{u_1} = \vect{\nabla}\psi_1 \cdot \vect{\mathcal{V}} = 0,
    \end{split} \quad \;\;\text{in} \quad \Omega,
\end{equation}

In order to make the boundary integral in Eq. \ref{appendix_perturbation_eq} zero, we require $\left(\vect{j_u} - \vect{\psi^T} \left(\vect{F_u}\cdot \vect{n}\right) \right)\vect{\delta u} = 0$ on $\Gamma$. Due to the prescribed boundary conditions in Eqs. \ref{u0_eq_boundary} and \ref{u1_eq_boundary}, $\delta u_0 = 0$ on $\Gamma_0$ and $\delta u_1 = 0$ on $\Gamma_1$. Thus, the adjoint solution needs to satisfy the following boundary conditions:
\begin{equation}
\label{adjoint_bc_prev}
    \begin{split}
      &j_{u_0} - \vect{\psi^T} \left(\vect{F}_{u_0}\cdot \vect{n}\right) = 0 \quad \text{on} \quad \Gamma\setminus \Gamma_0, \\  
      &j_{u_1} - \vect{\psi^T} \left(\vect{F}_{u_1}\cdot \vect{n}\right) = 0 \quad \text{on} \quad \Gamma\setminus \Gamma_1.
    \end{split}
\end{equation}
Substituting Eqs. \ref{FS_substituted} and \ref{appendix_functional} in \ref{adjoint_bc_prev} and noting that $\Gamma\setminus\Gamma_1 = \Gamma_{\text{right}} \cup \Gamma_{\text{bottom}}$, we get
\begin{equation}
    \begin{split}
        & \psi_0 \left(\vect{\beta}\cdot\vect{n}\right) = 0 \quad \text{on} \quad \Gamma\setminus \Gamma_0, \\
        & \psi_1 \left(\vect{\mathcal{V}}\cdot\vect{n}\right) = \begin{cases}
            2u_1, & \Gamma_{\text{right}}, \\
            0, & \Gamma_{\text{bottom}},
        \end{cases} 
    \end{split}
\end{equation}
giving Eqs. \ref{eq_psi0_boundary} and \ref{eq_psi1_boundary}.